\documentclass[journal]{IEEEtran}
\usepackage{mathtools}
\usepackage[hidelinks]{hyperref}
\usepackage{algorithmic}
\usepackage[ruled,vlined]{algorithm2e}       
\usepackage{amsthm,amsmath,amssymb}
\allowdisplaybreaks[4]
\usepackage{color}
\usepackage{cite}

\newtheorem{lemma}{Lemma}
\newtheorem{theorem}{Theorem}
\newtheorem{assumption}{Assumption}
\newtheorem{remark}{Remark}
\newtheorem{definition}{Definition}

\newtheorem{corollary}{Corollary}

\ifCLASSINFOpdf
\else
\fi

\begin{document}
%
\title{Constrained Optimization with \\ Decision-Dependent Distributions}

\author{
Zifan Wang, Changxin Liu, Thomas Parisini, Michael M. Zavlanos, and Karl H. Johansson
\thanks{This work was supported in part by Swedish Research Council Distinguished Professor Grant 2017-01078, Knut and Alice Wallenberg Foundation, Wallenberg Scholar Grant, the Swedish Strategic Research Foundation CLAS Grant RIT17-0046, AFOSR under award \#FA9550-19-1-0169, NSF under award CNS-1932011, an NSERC Postdoctoral Fellowship, the Digital Futures Scholar-in-Residence Program, the European Union's Horizon 2020 research and innovation programme under grant agreement no. 739551 (KIOS CoE), and the Italian Ministry for Research in the framework of the 2017 Program for Research Projects of National Interest (PRIN), Grant no. 2017YKXYXJ. (Corresponding author: Changxin Liu.)}
\thanks{Zifan Wang, Changxin Liu and Karl H. Johansson are with Division of Decision and Control Systems, School of Electrical Engineering and Computer Science, KTH Royal Institute of Technology, and also with Digital Futures, SE-10044 Stockholm, Sweden. Email: \{zifanw,changxin,kallej\}@kth.se.}
\thanks{Thomas Parisini is with the Department of Electrical and Electronic Engineering, Imperial College London, London SW7 2AZ, UK, with the Department of Engineering and Architecture, University of Trieste, 34127 Trieste, Italy, and also with the KIOS Research and Innovation Center of Excellence, University of Cyprus, CY-1678 Nicosia, Cyprus. Email: t.parisini@imperial.ac.uk}
\thanks{Michael M. Zavlanos is with the Department of Mechanical Engineering and Materials Science, Duke University, Durham, NC, USA. Email: michael.zavlanos@duke.edu}
}

\markboth{Journal of \LaTeX\ Class Files,~Vol.~14, No.~8, August~2015}%
{Shell \MakeLowercase{\textit{et al.}}:  Constrained Optimization with Decision-Dependent Distributions}
%



\maketitle

\begin{abstract}
In this paper we deal with stochastic optimization problems where the data distributions change in response to the decision variables.
Traditionally, the study of optimization problems with decision-dependent distributions has assumed either the absence of constraints or fixed constraints.
This work considers a more general setting where the constraints can also dynamically adjust in response to changes in the decision variables.
Specifically, we consider linear constraints and analyze the effect of decision-dependent distributions in both the objective function and constraints. 
Firstly, we establish a sufficient condition for the existence of a constrained equilibrium point, at which the distributions remain invariant under retraining.
Morevoer, we propose and analyze two algorithms: repeated constrained optimization and repeated dual ascent.
For each algorithm, we provide sufficient conditions for convergence to the constrained equilibrium point.
Furthermore, we explore the relationship between the equilibrium point and the optimal point for the constrained decision-dependent optimization problem.
Notably, our results encompass previous findings as special cases when the constraints remain fixed.
To show the effectiveness of our theoretical analysis, we provide numerical experiments on both a market problem and a dynamic pricing problem for parking based on real-world data.
\end{abstract}


%
\IEEEpeerreviewmaketitle

\section{Introduction}

Machine learning algorithms typically use historical data under the assumption that these data adequately reflect future system behavior. 
However, in many settings, decisions made by machine learning algorithms will influence behavior of the system and lead to changes to the data distribution. 
For example, in dynamic pricing problems for parking management \cite{pierce2018sfpark,dowling2019modeling}, the users may adjust their decisions on whether and for how long to park, causing modifications in the distributions of future data.
This phenomenon appears in numerous applications, such as online labor markets \cite{horton2010online}, predictive policing \cite{lum2016predict}, and vehicle sharing markets \cite{banerjee2015pricing,bianchin2021online,narang2022learning}, where strategic users react to changes of decisions and make adjustments accordingly.

To account for distribution shifts in the data in response to decisions influenced by machine learning models, new frameworks under the name of 
\textit{performative prediction} or \textit{optimization with decision-dependent distributions} have recently been developed \cite{perdomo2020performative}. 
Specifically, the term ``performative prediction" was coined in \cite{perdomo2020performative} to illustrate the phenomenon that predictive models can trigger reactions that influence the outcome they aim to predict. 
This framework also incorporates a theoretical analysis that relies on defining a distribution map that relates decisions to data distributions and associating this map  with a sensitivity property that employs the Wasserstein distance metric.

Many optimization problems with decision-dependent distributions also involve constraints that are themselves also influenced by the decisions. 
For example, in dynamic pricing for parking management problems, the total parking time, as a strategic feature in response to the change of the price (decision), makes the constraints decision-dependent if the total parking time plays a role in defining the constraints.
Motivated by such scenarios, in this paper we investigate constrained optimization with decision-dependent distributions affecting both the objective function and constraints. 
Specifically, we consider linear constraints that depend on the decision variables and introduce two notions of solution points: the constrained equilibrium point and constrained optimal point. Assuming that the distribution maps in both the objective function and constraints are insensitive, we show that the constrained equilibrium point is unique and repeated constrained minimization (RCM) converges to this unique equilibrium point. Notably, when the distributions in the constraints are independent of the decisions, the sufficient condition for convergence of RCM reduces to that for repeated risk minimization in [8]. As RCM requires the exact solution of a constrained optimization problem, we relax this requirement and propose the repeated dual ascent (RDA) algorithm. We show that convergence of RDA to the constrained equilibrium point is guaranteed under more restrictive conditions compared to RCM. Additionally, we present convergence analysis for both RCM and RDA when only the constraints are decision-dependent and the distribution in the objective function remains fixed.
Finally, we analyze the relation between the constrained equilibrium points and the constrained optimal points. Under mild conditions, we show that the distance between these points can be bounded. This result also encompasses the case when the constraints are fixed, as explored in [8]. To validate our theoretical analysis, we present numerical experiments on a market problem and a dynamic pricing problem for parking management that utilizes an open-sourced dataset.

To the best of our knowledge, constrained optimization with decision-dependent constraints has not been explored in the literature. 
Most closely related to our study is the vast literature on optimization with decision-dependent distributions that exclusively occur in the objective function \cite{perdomo2020performative,narang2022learning,mendler2020stochastic,miller2021outside,drusvyatskiy2023stochastic,ray2022decision,lu2023bilevel,zhao2022optimizing,li2022multi,li2022state,wood2023online}. In this context, there are no constraints or the constraints remain fixed throughout the optimization process. 
Among these works, the authors in \cite{perdomo2020performative} first provide a comprehensive framework for handling decision-dependent problems. Specifically, they define the equilibrium and optimal points, and show that repeated risk minimization and repeated gradient descent converge to the equilibrium point when the distribution map is less sensitive. 
Moreover, under certain conditions, \cite{perdomo2020performative} shows that the distance between equilibrium and optimal points can be bounded.
In the case that only stochastic gradient is available, the work in \cite{mendler2020stochastic} studies stochastic optimization problems with decision-dependent distributions and proposes two variants of the stochastic gradient method that converge to the equilibrium point.
The follow-up work in \cite{miller2021outside} focuses on directly optimizing the performative risk  and  proposes algorithms to find the optimal point.
Subsequent works have extended this framework to different scenarios. For example, \cite{narang2022learning} investigates multi-agent games with decision-dependent distributions and analyzes the convergence to the Nash equilibrium in strongly monotone games. 
Similarly, \cite{wang2023network} defines performative equilibrium equilibrium in multi-agent networked games and analyzes convergence to this equilibrium.
However, common in all the above works is that there are no constraints or the constraints are fixed.
The consideration of decision-dependent constraints introduces additional complexity to the optimization process, as the constraints now adapt to the decisions made, necessitating more sophisticated techniques to address the challenges posed by these adapting constraints.

Related to our work is also the literature on optimization with time-varying constraints \cite{subotic2021quantitative,cao2018online,yi2020distributed1,fazlyab2017prediction}. 
The reason is that when the decisions change during the optimization process, the constraints become effectively time-varying and impact the feasible solution space and further the optimal solution.
The work in \cite{subotic2021quantitative} derives a quantitative bound on the rate of change of the optimizer when the objective function and constraints are perturbed concurrently. This bound is then used to quantify the tracking performance for continuous-time algorithms in time-varying constrained optimization.
In the context of optimization with decision-dependent distributions,  all existing works assume static constraints, with the exception of \cite{wood2021online}.
In \cite{wood2021online}, the authors consider online optimization with decision-dependent distributions where the objective function and constraints vary over time, leading to time-varying equilibrium points. 
They propose a projected gradient descent method and demonstrate its ability to track the trajectory of the evolving equilibrium points.
However, all works discussed before, including \cite{wood2021online}, do not leverage the crucial information that the perturbations in the objective function and/or constraints are induced by the changes in the decisions. 
As a result, the techniques presented in these works cannot be used to analyze optimization problems with decision-dependent constraints as the ones considered here. 

The rest of this paper is organized as follows. In Section~\ref{sec:problem}, we formulate our problem and provide some preliminary results. In Section~\ref{sec:alg}, we propose two algorithms and analyze their convergence to the constrained equilibrium point. Section~\ref{sec:relation} analyzes the relation between the constrained equilibrium and optimal points. Section~\ref{sec:experiment} provides the numerical validation of the proposed algorithms for a market problem and a dynamic pricing problem for parking management. Finally, concluding remarks are given  in Section~\ref{sec:conclusion}.

\section{Problem Setup}\label{sec:problem}
In this section, we formulate the problem and provide some preliminary results. 
Throughout the paper, we let $\mathbb{R}^n$ denote $n$-dimensional Euclidean space with inner product $\langle \cdot, \cdot \rangle$ and 2-norm $\left\|\cdot \right\|$. We equip $\mathbb{R}^n$ with $\sigma$-algebra. For a matrix $A\in \mathbb{R}^{m\times n}$, we denote by $\left\| A\right\|_2$ the maximum singular value of the matrix $A$.
\subsection{Main Definitions}


Consider the following linearly constrained optimization problem:
\begin{align}\label{eq:problem}
    \min_x &\mathop{\mathbb{E}}_{z\sim \mathcal{D}(x)}[l(x,z)] \nonumber \\
    {\rm{s.t.}} & \quad Gx \leq \mathop{\mathbb{E}}_{w\sim \mathcal{D}_g(x)}[w],
\end{align}
where $x\in\mathbb{R}^n$ is the decision variable, $G\in \mathbb{R}^{d_w\times n}$ is a matrix associated with the linear constraint, and $z\in \mathbb{R}^{d_z}$ and $w\in \mathbb{R}^{d_w}$ are the random variables in the objective function and constraint, respectively. 
We assume that the loss function $l: \mathbb{R}^{n} \times \mathbb{R}^{d_z} \rightarrow \mathbb{R}$ is twice continuously differentiable in $(x,z)$.
The random variable $z$ follows the  distribution $\mathcal{D}(x)$ that depends on the decision variable $x$, where $\mathcal{D}$ is the distribution map from $\mathbb{R}^n$ to the space of distributions. 
%
We assume that suitable Borel measurability conditions hold so that the expected value operators $\mathop{\mathbb{E}}\limits_{z\sim \mathcal{D}(x)}[\cdot]$ and $\mathop{\mathbb{E}}\limits_{z\sim \mathcal{D}_g(x)}[\cdot]$ are well-defined, hence implying that the minimization operations are well-defined.
In addition, we assume that the constraints are decision-dependent, i.e., the random variable $w$ depends on the decision variable and follows the distribution $\mathcal{D}_g(x)$, where $\mathcal{D}_g$ is the distribution map from $\mathbb{R}^n$ to the space of distributions. 
Let $z$ and $w$ take values in the metric spaces $\mathcal{Z}$ and $\mathcal{W}$, respectively.
%
We define the solution points below.
\begin{definition}(Constrained Optimal Point).
A vector $x_o$ is called a constrained optimal point if it satisfies
\begin{align*}
    x_o =  \mathop{\rm{arg \; min}}_{x } & \; \mathop{\mathbb{E}}_{z\sim \mathcal{D}(x)}[l(x,z)] \nonumber \\
    {\rm{s.t.}} & \quad  Gx\leq \mathop{\mathbb{E}}_{w\sim \mathcal{D}_g(x)}[w] .
\end{align*}
\end{definition}
A constrained optimal point solves the problem \eqref{eq:problem}, but is usually hard to compute since the distribution maps $\mathcal{D}$ and $\mathcal{D}_g$ are generally unknown to the decision maker. Therefore, we introduce an alternative solution point below.

\begin{definition}(Constrained Equilibrium Point).
\label{def:CSP}
A vector $x_s$ is called a constrained equilibrium point if it satisfies
\begin{align*}
    x_s =  \mathop{\rm{arg \; min}}_{x}  & \mathop{\mathbb{E}}_{z\sim \mathcal{D}(x_s)}[l(x,z)] \nonumber \\
     {\rm{s.t.}}& \quad Gx\leq \mathop{\mathbb{E}}_{w\sim \mathcal{D}_g(x_s)}[w] .
\end{align*}
\end{definition}



The equilibrium point defined in Definition \ref{def:CSP} denotes a class of points, not necessarily optimal for \eqref{eq:problem}, that solve the constrained optimization problem using the distributions they induce. It is also referred to as a {stable} point in \cite{perdomo2020performative}.

In this work, we make the assumption that a constrained optimal point exists for \eqref{eq:problem}, which is a commonly accepted assumption as discussed in previous studies \cite{perdomo2020performative}. Furthermore, our analysis of two optimization methods will demonstrate the existence of equilibrium points under certain conditions on the loss function and the distribution map.


Our goal in this paper is to solve problem \eqref{eq:problem} by first determining the existence and uniqueness of constrained equilibrium points and next developing
algorithms that converge to the constrained equilibrium points. Moreover, we aim to explore the relation between the constrained equilibrium points and the optimal points.

\subsection{Key Assumptions and Preliminary Results}

The distribution maps in \eqref{eq:problem} play a crucial role in determining how the distributions respond to changes in decisions. 
In order to upper bound the rate of change at which the distributions change, we impose a Lipschitz condition on the distribution map, known as $\epsilon$-sensitivity.
\begin{definition}($\epsilon$-sensitivity).
We say that a distribution map $\mathcal{D}(\cdot)$ is $\epsilon$-sensitive if for all $x,x'$, we have 
\begin{align}
    W_1(\mathcal{D}(x),\mathcal{D}(x'))\leq \epsilon \left\|x-x'\right\|,
\end{align}
where $W_1$ denotes the earth mover's distance.
\end{definition}
The earth mover's distance, also known as $1$-Wasserstein distance, is a natural metric for quantifying the dissimilarity between probability distributions. According to Kantorovich-Rubinstein duality theorem \cite{kantorovich1958space}, the earth mover's distance of two probability distributions $P$ and $Q$ can be expressed as 
\begin{align*}
    W_1(P,Q)= \sup_{f \sim {\rm{Lip}_1}} \{ \mathbb{E}_{X\sim P}[f(X)] - \mathbb{E}_{Y\sim Q} [f(Q)]\},
\end{align*}
where ${\rm{Lip}_1}$ represents the set of $1$-Lipschitz continuous functions.
We make the following assumptions on the cost function and the distribution maps.
\begin{assumption}\label{assump:strong_convex}
The loss function $l(x,z)$ is $\gamma$-strongly convex in $x$ for every $z \in \mathcal{Z}$, i.e., 
\begin{align*}
    l(x,z) \geq l(x',z) + \langle \nabla_x l(x',z),x-x' \rangle + \frac{\gamma}{2}\left\|x-x'\right\|^2,
\end{align*}    
for all $x,x' \in \mathbb{R}^n$.
\end{assumption}

\begin{assumption}\label{assump:smooth}
$\nabla_x l(x,z)$ is $\beta_z$-Lipschitz continuous in $z$ for every $x\in \mathbb{R}^n$, i.e.,
\begin{align*}
    &\left\| \nabla_x l(x,z) - \nabla_x l(x,z')\right\| \leq \beta_z \left\| z-z'\right\|,
\end{align*}
for all $z,z'\in \mathcal{Z}$,
and $\nabla_x l(x,z)$ is $\beta_x$-Lipschitz continuous in $x$ for every $z\in \mathcal{Z}$, i.e.,
\begin{align*}
    &\left\| \nabla_x l(x,z) - \nabla_x l(y,z) \right\| \leq \beta_x \left\| x-y\right\|,
\end{align*}
for all $x,y\in \mathbb{R}^n$.
\end{assumption}

\begin{assumption} \label{assump:distribution_map}
The distribution maps $\mathcal{D}$ and $\mathcal{D}_g$ are $\epsilon$- and $\epsilon_g$-sensitive, respectively, namely,
\begin{align*}
    & W_1(\mathcal{D}(x),\mathcal{D}(y))\leq \epsilon \left\|x-y \right\|, \nonumber \\
    & W_1(\mathcal{D}_g(x),\mathcal{D}_g(y))\leq \epsilon_g \left\|x-y \right\|,
\end{align*}
for all $x,y \in \mathbb{R}^n$.
\end{assumption}
It is worth noting that Assumptions \ref{assump:strong_convex}--\ref{assump:distribution_map} are quite standard in the literature, e.g., \cite{perdomo2020performative,miller2021outside,mendler2020stochastic}. Moreover, we introduce the following assumption on the constraints.

\begin{assumption}\label{assumption:constraint}
    The matrix $G$ has full row rank.
\end{assumption}

Assumption \ref{assumption:constraint} is a technical condition that ensures the strong concavity of the Lagrange dual function, and the uniqueness of the optimal dual variable \cite{bazaraa2013nonlinear}. These two properties facilitate the contraction analysis of the proposed dual ascent algorithm. This assumption is commonly used in the literature regarding resource allocation \cite{wu2021new} and primal-dual optimization \cite{qu2018exponential}. 
Relaxing the rank condition would introduce additional complexity to the theoretical analysis, making it more challenging to establish rigorous results.
It is interesting to examine whether it can be relaxed to other constraint qualifications, e.g., linear independence constraint qualification (LICQ) or Slater condition; we leave it for future research.

In what follows, we present two important lemmas that will be used in the subsequent analysis. 
\begin{lemma}\cite{perdomo2020performative}\label{lemma:PP}
Suppose that the loss $l(x,z)$ is $\gamma$-strongly convex in $x$ for every $z\in \mathcal{Z}$, $\nabla_x l(x,z)$ is $\beta_z$-Lipschitz continuous in $z$, and the distribution map $\mathcal{D}(\cdot)$ is $\epsilon$-sensitive. Define $F(x')=\mathop{\rm{arg \; min}}\limits_{x \in \mathcal{X}} \mathop{\mathbb{E}}\limits_{z\sim \mathcal{D}(x')}[l(x,z)]$, where $\mathcal{X}$ is a convex set. For all $x,x' \in \mathcal{X}$, we have
\begin{align}
    \left\| F(x)-F(x')\right\| \leq \frac{\epsilon \beta_z}{\gamma} \left\| x - x' \right\|.
\end{align}
\end{lemma}
Lemma~\ref{lemma:PP} shows that the variation of the optimal solution across two distributions is bounded by the distance between the corresponding points that generate these distributions. However, Lemma~\ref{lemma:PP} is applicable only when the constraint is fixed. The following lemma explores the sensitivity of the optimal value with respect to perturbations of the constraint.

\begin{lemma}\cite{subotic2021quantitative}\label{lemma:sensitivity}
Consider the convex optimization problem: 
\begin{align}\label{eq:problem_lemma}
    \min_{x\in\mathcal{X}} & \quad f(x) \nonumber \\
     \mathrm{s.t.}  & \quad Gx \leq v, 
\end{align}
where $f$ is $\gamma$-strongly convex, $\beta$-smooth and twice continuously differentiable, $G \in \mathbb{R}^{d_w \times n}$ has full row rank, and $\mathcal{X}$ is a convex set. 
Suppose that $GG^{\rm{T}} \geq \underline{w}^2 I$ for some $\underline{w}>0$.
Then, there exists a continuous map $x^{*}(v)$ such that $x^{*}(v)$ is the optimum of the problem \eqref{eq:problem_lemma} and satisfies $ \left\| \nabla_v x^{*}(v) \right\| \leq l_{x^{*}} = \sqrt{\frac{\beta}{\gamma}}\frac{1}{\underline{w}}$.
\end{lemma}


\section{Algorithms that Converge to Equilibrium Points}\label{sec:alg}
This section focuses on the analysis of constrained equilibrium points.
Specifically, we discuss the existence and uniqueness of constrained equilibrium points. Moreover, we introduce two algorithms: RCM and RDA, and present sufficient conditions required for each algorithm to converge to the equilibrium point.

\subsection{Repeated Constrained Minimization}
We begin by analyzing RCM. We establish a sufficient condition for RCM to converge, which also ensures the uniqueness of the constrained equilibrium point.

\begin{algorithm}[t]
\caption{Repeated Constrained Minimization}
\begin{algorithmic}[1]
\STATE  \textbf{Input}: Initial variable $x_0$.
    \FOR {$t=0,1,2,\ldots$}
        \STATE $x_{t+1} =  \mathop{\rm{arg \; min}}_{x}  f_{x_t}(x)  \quad  {\rm{s.t.}} \quad Gx \leq \xi(x_t)$;
    \ENDFOR
\STATE  \textbf{Output}: Sequence $\{ x_t\}$
\end{algorithmic}\label{alg:RCM}
\end{algorithm}

Starting with an initial decision variable, RCM repeatedly solves a constrained optimization problem with the distribution induced by the previous decision. The RCM algorithm is presented in Algorithm~\ref{alg:RCM}. Specifically, at iteration $t+1$, we perform the following update:
\begin{align}\label{eq:RCM}
    x_{t+1} = & \mathop{\rm{arg \; min}}_{x}  f_{x_t}(x)  \quad  {\rm{s.t.}} \quad Gx \leq \xi(x_t),
\end{align}
where $f_{x_t}(x) :=\mathop{\mathbb{E}}\limits_{z\sim \mathcal{D}(x_t)}[l(x,z)] $ and $\xi(x_t):=\mathop{\mathbb{E}}\limits_{w\sim \mathcal{D}_g(x_t)}[w]$.
The convergence of RCM is closely related to the existence and uniqueness of the constrained equilibrium point. The main result is presented in the following theorem. The proof can be found in Appendix~\ref{app:RCM}.
\begin{theorem}\label{thm:RCM}
Suppose Assumptions~\ref{assump:strong_convex}--\ref{assumption:constraint} hold and let $L_{x^{*}} :=\sqrt{\frac{\beta_x}{\gamma  \lambda_{\min}(GG^{\rm{T}}) }}$. 
If
\begin{align}\label{eq:RCM:suff}
    \frac{\epsilon \beta_z}{\gamma}+ L_{x^{*}} \epsilon_g <1,
\end{align}
then problem \eqref{eq:problem} admits a unique equilibrium point $x_s$ and RCM converges  to $x_s$ at a linear rate.
\end{theorem}

Theorem~\ref{thm:RCM} demonstrates that RCM converges to the unique equilibrium point at a linear rate when the parameters $\epsilon$ and $\epsilon_g$ are small enough, i.e., the sensitivity of the optimization problem is small.
When the constraints are fixed ($\epsilon_g=0$), RCM is equivalent to repeated risk minimization (RRM) discussed in \cite{perdomo2020performative}. In this case, the sufficient condition in \eqref{eq:RCM:suff} required for RCM reduces to $\frac{\epsilon \beta_z}{\gamma}<1$, which aligns with the result in \cite{perdomo2020performative}. Therefore, RCM subsumes RRM in \cite{perdomo2020performative} as a special case.

\begin{figure}[t]
\begin{center}
\centerline{\includegraphics[width=0.7\columnwidth]{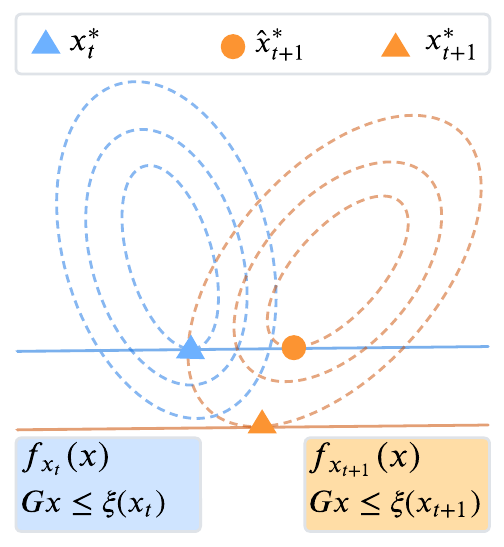}}
\caption{Illustration of an RCM iteration. The decision-dependent constraint shifts the location of the optimum.}
\label{RCM:illustration}
\end{center}
\vskip -0.2in
\end{figure}

When the constraint changes as a function of decisions, the dynamics of RCM become more intricate, as the solution  at each iteration may be influenced by the variation in the constraints. 
Fig.~\ref{RCM:illustration} illustrates the dynamics of RCM.
Under the constraint $Gx\leq \xi(x_t)$, the points $x_t^{*}$ and $\hat{x}_{t+1}^{*}$ represent minima of the objective functions $f_{x_t}(x)$ and $f_{x_{t+1}}(x)$, respectively. However, since the constraint shifts to $Gx\leq \xi(x_{t+1})$, the minimum at iteration $t+1$ is $x_{t+1}^*$ instead of $\hat{x}_{t+1}^{*}$. This constraint shift, which is usually hard to predict, makes the analysis of RCM more complex as shown in the proof.


When the objective function is fixed ($\epsilon=0$), implying that the distribution in the objective function is decision-independent, the convergence result is as follows. The proof follows directly from Theorem~\ref{thm:RCM}.
\begin{corollary}\label{corollary:RCM}
Suppose $\epsilon=0$ and Assumptions~\ref{assump:strong_convex}--\ref{assumption:constraint} hold, and let $L_{x^{*}} :=\sqrt{\frac{\beta_x}{\gamma  \lambda_{\min}(GG^{\rm{T}}) }}$.
If $ L_{x^{*}} \epsilon_g <1$, then problem \eqref{eq:problem} admits a unique equilibrium point $x_s$ and RCM converges  to $x_s$ at a linear rate.
\end{corollary}

In the case that the distribution in the objective function remains fixed, the convergence of RCM is thus ensured under the condition $ L_{x^{*}} \epsilon_g <1$.
The parameter $L_{x^{*}}$ can be interpreted as an intrinsic property of the system that characterizes the sensitivity to perturbations in the constraints. 
The parameter $\epsilon_g$ quantifies the extent to which the constraints are sensitive to changes in the decisions. 
For convergence, it is essential that both $L_{x^{*}}$ and $\epsilon_g$ are sufficiently small. 
In other words, it is important that the system exhibits low sensitivity to perturbations in the constraints, while simultaneously ensuring that the constraints are also minimally affected by changes in decisions. 
This ensures that the sequence of optimal solutions remains insensitive to the changes in the decision-dependent distributions.

When $\epsilon=0$, the sufficient condition $ L_{x^{*}} \epsilon_g <1$ for RCM to converge is in fact tight. This can be illustrated in the following simple example.

\noindent  \textbf{Example.} 
Consider the minimization problem with the loss function $x^2$ and the constraint $x\geq \xi(x_t)$, where $x\in \mathbb{R}$ is the decision variable, and $\xi(x_t) = \theta x_t$ with $\theta>0$. In this problem, it is easy to verify that $\epsilon=0$, $\epsilon_g = \theta$, $\gamma=\beta_x =2$, $G=1$, and $L_{x^{*}}=1$. If we run RCM from a positive initial point $x_0>0$, it is easy to verify that $x_t = \theta^t x_0$. The convergence of RCM requires the condition $\theta<1$, which is equivalent to $L_{x^{*}} \epsilon_g <1$.

\begin{remark}
Although the sufficient condition \eqref{eq:RCM:suff} seems to be conservative, it is in fact tight in two specific scenarios: when $\epsilon=0$ and when $\epsilon_g=0$, as evidenced by Proposition 3.6 in \cite{perdomo2020performative} and the example provided above, respectively. Since this condition is applicable to a wide range of distribution maps, relaxing it when $\epsilon$ and $\epsilon_g$ are both non-zero poses a considerable challenge, which is left as our future work.
\end{remark}

\subsection{Repeated Dual Ascent}
Theorem~\ref{thm:RCM} analyzes the convergence of RCM when the sensitivity parameters $\epsilon,\epsilon_g$ are small enough. 
However,  implementing RCM requires the exact solution of a constrained optimization problem at each iteration, which might be computationally inefficient. To address this limitation, we  introduce a dual ascent algorithm RDA in this section.
Sufficient conditions are proposed that guarantee the convergence of this algorithm to the constrained equilibrium point.

Before presenting the algorithm, we define the Lagrangian function
\begin{align*}
    L_{x'}(x,\lambda):=f_{x'}(x) + \lambda^{\rm{T}}(Gx - \xi(x'))
\end{align*}
where $f_{x'}(x) := \mathop{\mathbb{E}}\limits_{z\sim \mathcal{D}(x')}[l(x,z)]$, and the dual function
\begin{align*}
    d_{x'}(\lambda):= \mathop{\min}\limits_{x} L_{x'}(x,\lambda).
\end{align*}
It is easy to verify that the loss function $f_{x'}(x)$ is $\gamma$-strongly convex and $\beta_x$-smooth in $x$ for every $x'\in \mathbb{R}^n$ by taking the expectation of the loss function $l(x,z)$ with respect to $z\sim\mathcal{D}(x')$.


Note that in this setting the dual function can be written as 
\begin{align}
    d_{x'}(\lambda) = - \tilde{f}_{x'}(-G^{\rm{T}}\lambda) - \lambda^{\rm{T}} \xi(x'),
\end{align}
where the conjugate function $\tilde{f}_{x'}$ is defined as $\tilde{f}_{x'}(y):= \mathop{\max}\limits_x  \langle y, x \rangle - f_{x'}(x)$. The gradient of the dual function is given by
\begin{align}
    \nabla d_{x'}(\lambda)=G \nabla \tilde{f}_{x'}(-G^{\rm{T}} \lambda ) - \xi(x').
\end{align}

We present the RDA as Algorithm~\ref{alg:RDA}. Specifically, at iteration $t$, the primal variable $x_t$ is obtained by minimizing the Lagrangian function $L_{x_{t-1}}(x,\lambda_t)$, i.e.,
\begin{align}
    &x_{t} =  \mathop{\rm{arg \; min}}_{x} L_{x_{t-1}}(x,\lambda_t).
\end{align}
Then, the dual variable is updated via gradient ascent:
\begin{align}\label{eq:RDA:dual:update}
    \lambda_{t+1} = [\lambda_t + \eta \nabla d_{x_t}(\lambda_t)]_{+}  = [\lambda_t + \eta (G\bar{y}_t -\xi(x_t) )]_{+},
\end{align}
where $\bar{y}_t = \mathop{\rm{arg \; min}}\limits_{x} L_{x_t}(x,\lambda_t)$ and $[a]_{+} = \max\{ a,0\}$. 
Note that the update direction in \eqref{eq:RDA:dual:update} is the positive gradient of the dual function $d_{x_t}$ instead of $d_{x_{t-1}}$. In order to compute the gradient of the new dual function $\nabla d_{x_t}(\lambda_t)$, we need to find the minima $\bar{y}_t$ by minimizing the loss function $L_{x_t}(x,\lambda_t)$.

\begin{algorithm}[t]
\caption{Repeated Dual Ascent}
\begin{algorithmic}[1]
\STATE  \textbf{Input}: Initial variables $\lambda_1$, $x_0$, step size $\eta$.
    \FOR {$t=1,2,\ldots$}
        \STATE $x_{t} =  \mathop{\rm{arg \; min}}\limits_{x} L_{x_{t-1}}(x,\lambda_t)$;
        \STATE $\bar{y}_t = \mathop{\rm{arg \; min}}\limits_{x} \{f_{x_t}(x) + \lambda_t^{\rm{T}} (Gx - \xi(x_t))\}$;
        \STATE $\lambda_{t+1}=[\lambda_t + \eta (G\bar{y}_t -\xi(x_t) )]_{+}$;
    \ENDFOR
\STATE  \textbf{Output}: Sequences $\{ x_t\}$, $\{ \lambda_t \}$
\end{algorithmic}\label{alg:RDA}
\end{algorithm}

Given that $f_{x'}(x)$ is $\gamma$-strongly convex and $\beta_x$-smooth in $x$ for every $x'\in \mathbb{R}^n$, 
it can be shown that the dual function $d_{x'}(\lambda)$ is strongly concave and smooth. This result is summarized in the following lemma and the proof can be found in the literature; see Proposition 3.1 in \cite{guigues2020strong} for the analysis of strong concavity and Lemma 3.2 in \cite{luo1993convergence} for the analysis of smoothness.
\begin{lemma}
Suppose Assumptions~\ref{assump:strong_convex}, \ref{assump:smooth}, and \ref{assumption:constraint} hold. There holds that, for every $x'\in\mathbb{R}^n$, $d_{x'}(\lambda)$ is $\gamma_d$-strongly concave, i.e.,
\begin{align*}
    d_{x'}(\lambda) \leq d_{x'}(\lambda') + \langle \lambda-\lambda', \nabla d_{x'}(\lambda')\rangle + \frac{\gamma_d}{2}\left\| \lambda- \lambda'\right\|^2,
\end{align*}
with $\gamma_d=\frac{\lambda_{\min}(G G^{\rm{T}})}{L}$, and $L_d$-smooth, i.e.,
\begin{align*}
    \left\| \nabla d_{x'}(\lambda) - \nabla d_{x'}(\lambda') \right\|\leq L_d \left\| \lambda - \lambda'\right\|,
\end{align*}
with $L_d = \frac{\left\| G \right\|_2^2 }{\gamma}$.
\end{lemma}

We say that the linear independence constraint qualification (LICQ) is satisfied if the gradients of all active constraints are linearly independent.
The following lemma discusses LICQ and the uniqueness of the optimal dual variable.
\begin{lemma}\label{lemma:LICQ} \cite{bazaraa2013nonlinear} 
Suppose Assumptions~\ref{assump:strong_convex} and \ref{assumption:constraint} hold. Then, the LICQ holds and the optimal dual variable is unique.
\end{lemma}
Now we are ready to present the convergence result of RDA. The proof can be found in Appendix~\ref{app:RDA}.

\begin{theorem}\label{thm:RDA}
Suppose Assumptions~\ref{assump:strong_convex}--\ref{assumption:constraint} hold, and denote by $\{(x_t,\lambda_t)\}$ the sequence generated by RDA.
If
\begin{align}\label{eq:RDA:suff1}
    \left(\frac{\epsilon \beta_z \left\| G \right\|_2}{\gamma} +\epsilon_g \right)\left(1+\frac{ \epsilon \beta_z \left\| G \right\|_2}{\gamma^2} + \frac{\left\| G \right\|_2^2}{\gamma^2} \right)< 2\gamma_d
\end{align}
and
\begin{align}\label{eq:RDA:suff2}
    \frac{\epsilon \beta_z}{\gamma}\left(1+  \frac{\epsilon \beta_z \left\| G \right\|_2}{\gamma^2} \right)<1
\end{align}
hold,
then there exist $s_2>s_1>0$ such that for all step size $\eta \in (s_1,s_2)$ we have
\begin{align}
    &\left\| \lambda_{t+1} - \lambda_s^{*} \right\|^2 + \alpha \left\| x_t - x_s \right\|^2 \nonumber \\
    & \leq \kappa^t (\left\| \lambda_{1} - \lambda_s^{*} \right\|^2 + \alpha \left\| x_0 - x_s \right\|^2),
\end{align}
where $\lambda_s^{*}$ is the optimal dual variable of the dual function $d_{x_s}(\lambda)$, and 
\begin{align*}
    &\rho_1 = \frac{\epsilon \beta_z}{\gamma}, \;\rho_2 = \rho_1 \left\| G \right\|_2 +\epsilon_g, \\
    & \rho_3 =\frac{\rho_1 \left\| G \right\|_2}{\gamma} + \frac{\left\| G \right\|_2^2}{\gamma^2},\; \rho_4 = \rho_1^2 +\frac{\rho_1 \left\| G \right\|_2}{\gamma},\\
    & a_1 = L_d^2 + L_d \rho_2 + \rho_2^2 \rho_3 + L_d \rho_2 \rho_3, \\
    & a_2 = 2\gamma_d -\rho_2(1+\rho_3),  \; a_3 = \alpha \rho_3, \\
    & s_1 = \frac{a_2 - \sqrt{a_2^2 - 4a_1 a_3}}{2a_1}, \\
    &b_1  = \rho_4(\rho_2^2 + L_d \rho_2), \; b_2 = \rho_2 \rho_4, \; b_3=\alpha(1-\rho_4), \\
    &s_2 = \frac{-b_2 + \sqrt{b_2^2 +4 b_1 b_3}}{2b_1},\; \alpha = (1-\sqrt{\bar{\alpha}})\frac{a_2^2}{4\rho_3 a_1}, \\
    & \sqrt{\bar{\alpha}} = \max\left\{0,1-\frac{L_d^2}{4\gamma_d(\rho_2 +L_d)}\right\} ,\\
    &\kappa = \max (\kappa_1,\kappa_2), \;
    \kappa_1 = \eta^2 a_1  - \eta a_2 + \alpha a_3<1, \\
    &\kappa_2 = \eta^2 b_1 +\eta b_2 + b_3<1.
\end{align*}
\end{theorem}

Theorem~\ref{thm:RDA} shows that RDA converges to the equilibrium point at a linear rate under conditions \eqref{eq:RDA:suff1} and \eqref{eq:RDA:suff2}.
Note that the sufficient condition in \eqref{eq:RDA:suff1} is more restrictive than the condition in \eqref{eq:RCM:suff} that guarantees the convergence of RCM to the equilibrium point, because
\begin{align*}
    &\frac{1}{2\gamma_d} \left(\frac{\epsilon \beta_z}{\gamma}\left\| G \right\|_2 +\epsilon_g \right)\left(1+\frac{ \epsilon \beta_z \left\| G \right\|_2}{\gamma^2} + \frac{\left\| G \right\|_2^2}{\gamma^2} \right) \\
    & \geq \frac{1}{2\gamma_d} \left(\frac{\epsilon \beta_z}{\gamma}\left\| G \right\|_2 +\epsilon_g \right) \frac{2\left\| G\right\|_2 }{\gamma} \\
    & = \frac{\epsilon \beta_z}{\gamma} \frac{ \left\|G\right\|_2^2 \beta_x }{\gamma \lambda_{\min}(G G^T)} + L_{x^{*}} \epsilon_g \sqrt{\frac{ \left\| G\right\|_2^2 \beta_x}  { \gamma \lambda_{\min}(G G^T) } }  \\
    &\geq \frac{\epsilon \beta_z}{\gamma} + L_{x^{*}} \epsilon_g,
\end{align*}
where in the first inequality we use $1+ \frac{\left\| G \right\|_2^2}{\gamma^2} \geq \frac{2\left\| G \right\|_2 }{\gamma}$ and in the equality we use the definitions $\gamma_d=\frac{\lambda_{\min}(G G^{\rm{T}})}{\beta_x}$ and $ L_{x^{*}} = \sqrt{\frac{\beta_x}{\gamma \lambda_{\min}(G G^T)} }$. The last inequality follows from the facts that the strong convexity parameter is not larger than the smooth parameter, i.e., $\gamma \leq \beta_x$ and $\lambda_{\min}(G G^{\rm{T}}) \leq \left\| G\right\|_2^2$.
Moreover, it follows that the condition \eqref{eq:RDA:suff1} together with Assumptions~\ref{assump:strong_convex}--\ref{assumption:constraint} ensures the uniqueness of the equilibrium point.

\begin{remark}
When the distribution in the constraints is fixed, one can use the repeated gradient descent (RGD) method with projection in \cite{perdomo2020performative}, which performs gradient descent in the primal space.
It can be verified that a sufficient condition for RGD to converge to the equilibrium point is $ \frac{\epsilon \beta_z}{\gamma}<1$. 
Interestingly, we observe that RDA, which performs gradient ascent in the dual space, requires more restrictive conditions compared to RGD due to the fact that one of the sufficient conditions \eqref{eq:RDA:suff2} for RDA to converge is more stringent than the condition $ \frac{\epsilon \beta_z}{\gamma}<1$.
\end{remark}

\begin{remark}
We note that RDA deviates from the traditional dual ascent algorithm by involving the solution of two minimization problems at each iteration, as opposed to the traditional algorithm that only requires one. The additional computation of $\bar{y}_t$ is implemented to ensure that the dual variable performs an update in the direction of the updated gradient $ \nabla d_{x_t}(\lambda_t)$. An alternative way to update the dual variable is to use $ \nabla d_{x_{t-1}}(\lambda_t)= G x_t -\xi(x_{t-1})$, which eliminates the need for computing $\bar{y}_t$. However, the theoretical analysis of the algorithm with this update becomes challenging and we leave it as a subject for future work. 
\end{remark}

When the objective function is fixed ($\epsilon=0$), implying that the distribution in the objective function is decision-independent, the convergence of RDA is presented in the following result. The proof can be found in Appendix~\ref{app:corollary2}.
\begin{corollary}\label{corollary:RDA}
Suppose $\epsilon =0$ and Assumptions~\ref{assump:strong_convex}--\ref{assumption:constraint} hold, and denote by $\{(x_t,\lambda_t)\}$ the sequence generated by RDA.
If
\begin{align}\label{eq:RDA:epsi0:suff}
    \epsilon_g \left(1 + \frac{\left\| G \right\|_2^2}{\gamma^2} \right)< 2\gamma_d
\end{align} 
and
\begin{align*}
    \eta <\frac{2\gamma_d - \epsilon_g (1 + \frac{\left\| G \right\|_2^2}{\gamma^2} )}{L_d^2 + \epsilon_g L_d + \epsilon_g^2 \frac{\left\| G \right\|_2^2}{\gamma^2} + \epsilon_g^2 L_d \frac{\left\| G \right\|_2^2}{\gamma^2} },
\end{align*} 
then there holds
\begin{align}\label{eq:corollary:lambda_converge}
    &\left\| \lambda_{t+1} - \lambda_s^{*} \right\|^2  
    \leq \kappa_3^t \left\| \lambda_{1} - \lambda_s^{*} \right\|^2
\end{align}
and
\begin{align}\label{eq:corollary:x_converge}
    \left\| x_{t+1} - x_s \right\|^2  
    \leq \kappa_3^{t-1} \frac{\left\| G\right\|_2^2 }{\gamma^2} \left\| \lambda_{1} - \lambda_s^{*} \right\|^2,
\end{align}
where $\kappa_3 = 1- 2\eta \gamma_d + \eta^2 L_d^2 + \epsilon_g \eta +\epsilon_g L_d \eta^2 + \frac{\left\| G \right\|_2^2}{\gamma^2} ( \epsilon_g^2 \eta^2 +\epsilon_g \eta + \epsilon_g L_d \eta^2)$.
\end{corollary}

Condition \eqref{eq:RDA:epsi0:suff} required for RDA to converge is more stringent than the condition $L_{x^{*}} \epsilon_g <1$ that ensures the convergence of RCM, because
\begin{align*}
    &\frac{1}{2 \gamma_d} \epsilon_g \left(1 + \frac{\left\| G \right\|_2^2}{\gamma^2} \right)  \geq  \frac{1}{2 \gamma_d} \epsilon_g \frac{ 2\left\| G \right\|_2 }{\gamma} \nonumber \\
    & \geq \epsilon_g L_{x^{*}} \sqrt{\frac{\beta_x \left\| G\right\|_2^2 }{\gamma \lambda_{\min}(G G^T)} }\geq \epsilon_g L_{x^{*}},
\end{align*}
where the last inequality holds since $\beta_x \geq \gamma$ and $\left\| G\right\|_2^2 \geq \lambda_{\min}(G G^T) $. This aligns with the case when $\epsilon>0$.

\begin{remark}
Theorem~\ref{thm:RDA} demonstrates that the step size $\eta$ cannot be arbitrarily small to guarantee convergence of RDA.
In contrast, the traditional dual ascent method allows for an arbitrarily small step size.
This requirement arises here since a small step size may hinder the dual variable from keeping pace with changes in the optimal dual variable which are induced by changes in distributions in the objective function. From a technical perspective, a small step size cannot guarantee that the term $\left\| \lambda_t - \lambda_s^{*} \right\|^2 $ will decrease with time. However, when the distribution in the objective function is fixed ($\epsilon=0$), this concern is eliminated, and the step size can be selected arbitrarily small, as shown in Corollary~\ref{corollary:RDA}.
\end{remark}


\section{Relating Equilibrium and Optimal Points}\label{sec:relation}
The results presented so far focus on the convergence to the equilibrium points. However, a natural question to consider is the relation between the constrained equilibrium and optimal points.
In this section, we quantify the distance between them.
To proceed, we introduce the following additional assumption.
\begin{assumption}\label{assump:Lz}
$l(x,z)$ is $L_z$-Lipschitz continuous in z for every $x$.
\end{assumption}
Assumption~\ref{assump:Lz} is common when analyzing the relation between equilibrium and optimal points, e.g., \cite{perdomo2020performative}.
Given Assumption~\ref{assump:Lz}, our next result shows that all the equilibrium and optimal points are close to each other.
The proof can be found in Appendix~\ref{app:optimality_stability}.
\begin{theorem}\label{thm:optimality_stability}
Suppose Assumptions~\ref{assump:strong_convex}--\ref{assump:distribution_map} and \ref{assump:Lz} hold. If $L_{x^{*}} \epsilon_g < 1$, then for every constrained equilibrium point $x_s$ and every optimal point $x_o$, we have 
\begin{align}\label{eq:thm3:bound}
    \left\| x_o - x_s \right\| \leq \frac{2(L_z \epsilon + \epsilon_g\sqrt{d_w} \left\| \lambda_s^{*} \right\|)}{\gamma (1-L_{x^{*}} \epsilon_g \sqrt{d_w})^2}.
\end{align}
\end{theorem}

Note that Theorem~\ref{thm:optimality_stability} remains valid also for cases when constrained equilibrium points and constrained optimal points are not necessarily unique. 
Theorem~\ref{thm:optimality_stability} shows that all equilibrium and optimal points are close to each other if the parameters $\epsilon$ and $\epsilon_g$ are small. Consequently, constrained optimal points can be approximated by constrained equilibrium points which can be derived using RCM and RDA when $\epsilon$ and $\epsilon_g$ meet the specified conditions. 

When $\epsilon_g =0$, the bound in \eqref{eq:thm3:bound} becomes equivalent to the result in \cite{perdomo2020performative}. Therefore, our result subsumes the result in \cite{perdomo2020performative} as a special case. This demonstrates the broader applicability  of our framework.

To end this section, we briefly discuss the challenges encountered in the analysis of Theorem~\ref{thm:optimality_stability}. 
When analyzing the distance between the equilibrium and optimal points, a standard technique is to use the strong convexity property, which is also used in our proof. However, this property cannot be applied in a straightforward way here since the constraint sets induced by the distributions at the constrained equilibrium and optimal points are not the same. 
The presence of inconsistent constraints necessitates us to define a new optimal point, subject to newly constructed constraints that ensure that all these points satisfy the same constraints. 
By introducing this intermediate point and conducting the sensitivity analysis, we are able to address the challenges arising from the presence of inconsistent constraints.
Specifically, we can establish an upper bound on the distance between the equilibrium and optimal points and provide valuable insights into the relation between the equilibrium and optimal solutions for the constrained problem considered here.

\section{Numerical Experiments}\label{sec:experiment}
In this section, we conduct experiments to illustrate the performance of RCM and RDA on market and dynamic pricing problems. 
\subsection{Market Problem}
We consider a market problem with the cost function $l(x,z)=-x_1 z_1 - x_2 z_2 $. Here, $x_i$ is the price of the $i$-th good and $z_i$ is the random demand of the $i$-th good, $i=1,2$. 
We assume that the demand of the first good satisfies
$z_1 = \zeta_1(x_1) - a_1 x_1 $, where the random variable $\zeta_1(x_1)$ follows the uniform distribution $U[\zeta_{L_1},\zeta_{R_1}+ \epsilon x_1]$. 
For the second good, the random demand is given by $z_2 = \zeta_2 - a_2 x_2$, where $\zeta_2$ follows a fixed uniform distribution $U[\zeta_{L_2}, \zeta_{R_2}]$.
Additionally, we introduce the cost of production $v_i$ for each good. We assume that the cost of the first good, $v_1$, follows a uniform distribution $U[\underline{v}_1, 1.2 \underline{v}_1 + \epsilon_g x_1]$, where a lower price leads to higher demands and, consequently, lower average cost.
The cost of the second good, $v_2$, follows a uniform distribution $U[ \underline{v}_2, 1.2 \underline{v}_2]$.
The objective is to maximize the gross sales, represented by $x_1 z_1 + x_2 z_2$. 
We consider the constraint $Gx\leq \mathop{\mathbb{E}}[\xi(x)]$, where $G=[-a_3, -a_4]$, $\xi(x)=-v_1(x_1)-v_2 -e_1$, which ensures that the selling price is greater than the weighted average cost plus $e_1$.

The decision-dependent constrained optimization problem is defined as 
\begin{align*}
    \min& \quad - x_1 ( \mathbb{E}[\zeta_1(x_1)] - a_1 x_1 ) - x_2 (\mathbb{E}[\zeta_2] - a_2 x_2) \\
     {\rm{s.t.}} & \quad a_3 x_1 + a_4 x_2 \geq e_1 + \mathbb{E}[v_1(x_1)] + \mathbb{E}[v_2].
\end{align*}
The distributions in both the objective function and constraints depend on the decisions. Next, we analyze the sensitivity of the distribution maps.
Without loss of generality, we assume $x_1'>x_1$.
By virtue of the definition of Wasserstein distance, we have
\begin{align*}
    &W_1( v_1(x_1), v_1(x_1')) \nonumber \\
    &=1-\frac{\eta_{R_1} +\epsilon x_1 - \zeta_{L_1} }{\zeta_{R_1}+\epsilon x_1' - \zeta_{L_1} } + 1 - \frac{\zeta_{R_1} + \epsilon x}{\zeta_{R_1} + \epsilon x' } \\
    &\leq \epsilon |x_1 - x_1'| \left(\frac{1}{\zeta_{R_1}+\epsilon x_1' - \zeta_{L_1} } + \frac{1}{ \zeta_{R_1} + \epsilon x' }\right) \\
    &\leq \epsilon \left(\frac{1}{\zeta_{R_1} - \zeta_{L_1}} + \frac{1}{\zeta_{R_1}}\right) |x_1 - x_1'|.
\end{align*}
Based on the last inequality, we select the sensitivity parameter of the distribution map in the objective function as $\bar{\epsilon} = \epsilon (\frac{1}{\zeta_{R_1} - \zeta_{L_1}} + \frac{1}{\zeta_{R_1}})$. Similarly, for the sensitivity parameter of the distribution map of the constraints, we set $\bar{\epsilon}_g = \epsilon_g (\frac{1}{0.2 \underline{v}_1} +\frac{1}{1.2 \underline{v}_1} )$.

\begin{figure}[t]
\begin{center}
\centerline{\includegraphics[width=0.9\columnwidth]{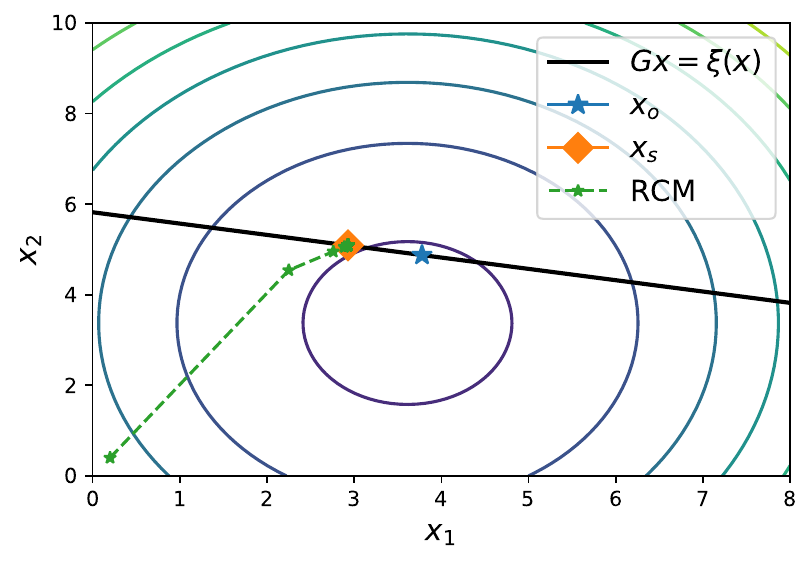}}
\caption{Trajectory of RCM when $\epsilon=0.7$, $\epsilon_g=0.7$ in the market problem. As expected, it conerges quickly to the equilibrium point $x_s$, which is close to the optimal point $x_o$.}
\label{RCM:Traj}
\end{center}
\vskip -0.2in
\end{figure}

\begin{figure}[t]
\begin{center}
\centerline{\includegraphics[width=0.9\columnwidth]{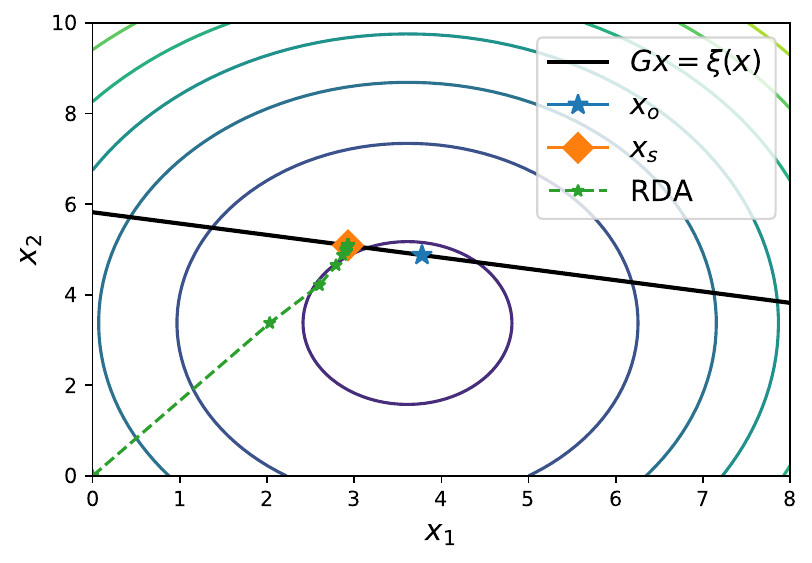}}
\caption{Trajectory of RDA when $\epsilon=0.7$, $\epsilon_g=0.7$ in the market problem. Slightly slower convergence compared to RCM in Fig.~\ref{RCM:Traj}.}
\label{RDA:Traj}
\end{center}
\vskip -0.2in
\end{figure}

\begin{figure}[t]
\begin{center}
\centerline{\includegraphics[width=0.9\columnwidth]{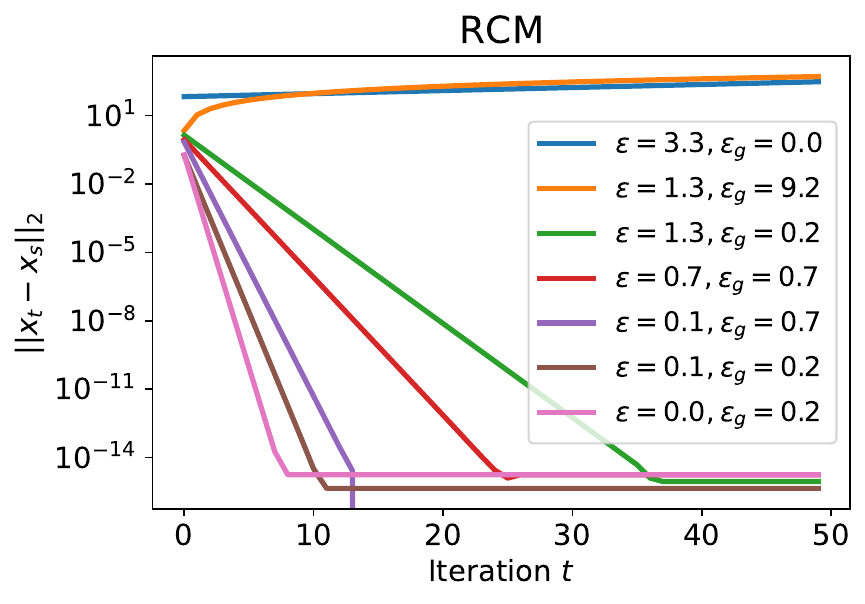}}
\caption{Convergence of RCM for different $\epsilon$, $\epsilon_g$ parameters in the market problem.
}
\label{RCM:market}
\end{center}
\vskip -0.2in
\end{figure}

\begin{figure}[t]
\begin{center}
\centerline{\includegraphics[width=0.9\columnwidth]{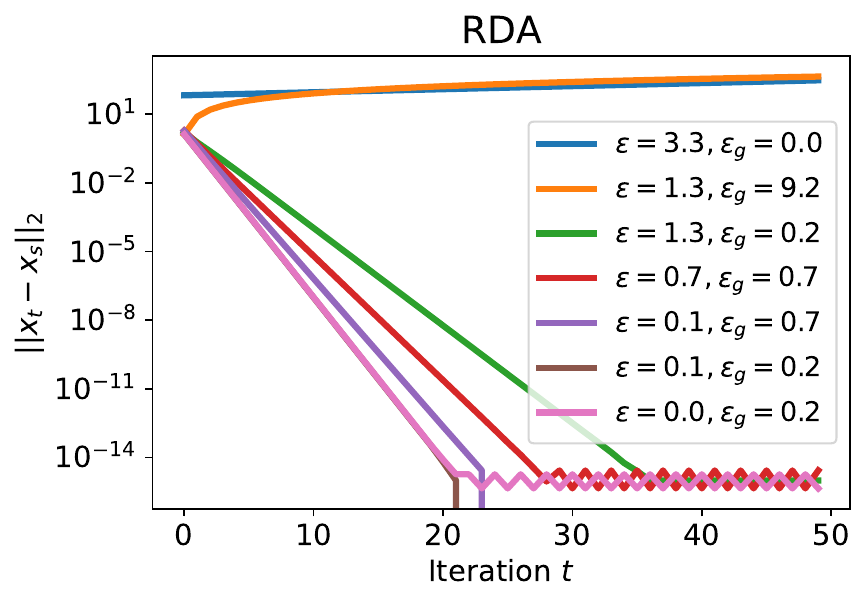}}
\caption{Convergence of RDA for different $\epsilon$, $\epsilon_g$ parameters in the market problem.}
\label{RDA:market}
\end{center}
\vskip -0.2in
\end{figure}

We consider the following parameter values for the market problem: $a_1=0.8$, $a_2=0.2$, $a_3=0.6$, $a_4 =1$, $\zeta_{L_1}=1$, $\zeta_{R_1}=5.5$, $\zeta_{L_2}=0.5$, $\zeta_{R_2}=2.2$, $e_1=1.2$, $\underline{v}_1=1.7$, $\underline{v}_2=2.5$. 
All the minimization problems are solved using the CVXOPT toolbox.
The simulation results for RCM and RDA are presented in Figs~\ref{RCM:Traj}--\ref{RDA:market}. In particular, Figs.~\ref{RCM:Traj} and \ref{RDA:Traj} show the evolution of the decision variable $x$ for RCM and RDA, respectively. 
In Figs.~\ref{RCM:market} and \ref{RDA:market}, we plot the performance of RCM and RDA for different values of the sensitivity parameters. 
We observe that when $\epsilon$ or $\epsilon_g$ is large, both methods fail to converge to the equilibrium point. 
When $\epsilon$ and $\epsilon_g$ are small, both methods converge rapidly within a few iterations. 
Moreover, when the sensitivity parameter $\epsilon$ is not that small, RCM performs similarly to RDA. However, when $\epsilon$ is small, RCM outperforms RDA, demonstrating faster convergence.

\subsection{Dynamic Pricing Problem}
In this subsection, we present the application of our algorithms to a dynamic pricing experiment using real-world parking data from SFpark in San Francisco \cite{sfpark_data}. The decision to use a personal vehicle for a trip is heavily influenced by parking availability, parking location, and price.
The primary goal of the SFpark pilot project was to make it easy to find a parking space. To this end,
SFpark implemented the following adjustments to hourly rates based on occupancy levels:
a) When the occupancy is between 80\% and 100\%, rates are increased by \$0.25;
b) When the occupancy is between 60\% and 80\%, no adjustment is made to the rates;
c) When the occupancy is between 30\% and 60\%, rates are decreased by \$0.25;
d) When the occupancy is below 30\%, rates are decreased by \$0.50.
SFpark aims to maintain occupancy $z_1$ between 60\% and 80\% and at the same time maximize the revenue $z_2(x+ \bar{x})$, where $z_2$ represents the total parking time and $x$ represents the difference between the parking price and the nominal price $\bar{x}=3$. 
We use a parameter $t>0$ to quantify the trade-off between these two objectives.  
%
The loss function for this dynamic pricing problem is defined as follows:
$$l(x,z_1,z_2)= (z_1 - 0.7)^2 - t z_2(x+ \bar{x}) + \frac{v}{2} \left\|x\right\|^2.$$
Here, $v>0$ is a regularization parameter.
When analyzing the individual response to a price adjustment $x$, we assume that each user updates their behavior using the best response strategy: $z_2' = z_2 - \epsilon x$. 
Here, $\epsilon$ is a positive constant that quantifies the sensitivity of the distribution map. It is easy to verify that the distribution map for $z_2$ is $\epsilon$-sensitive. 
\begin{figure}[t]
\begin{center}
\centerline{\includegraphics[width=0.9\columnwidth]{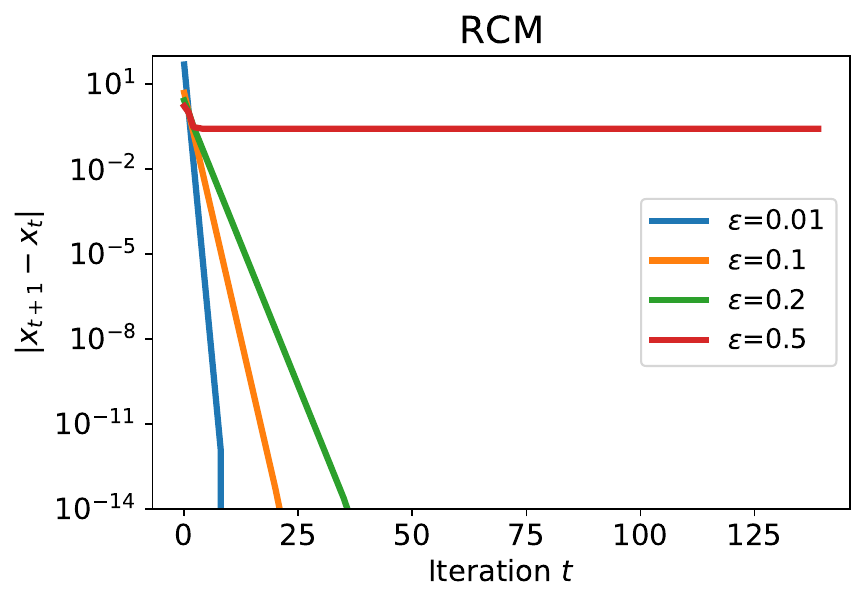}}
\caption{Convergence of RCM for different $\epsilon$-sensitivity parameters in the dynamic parking pricing problem.}
\label{RCM:pricing_simulation}
\end{center}
\vskip -0.2in
\end{figure}

\begin{figure}[t]
\begin{center}
\centerline{\includegraphics[width=0.9\columnwidth]{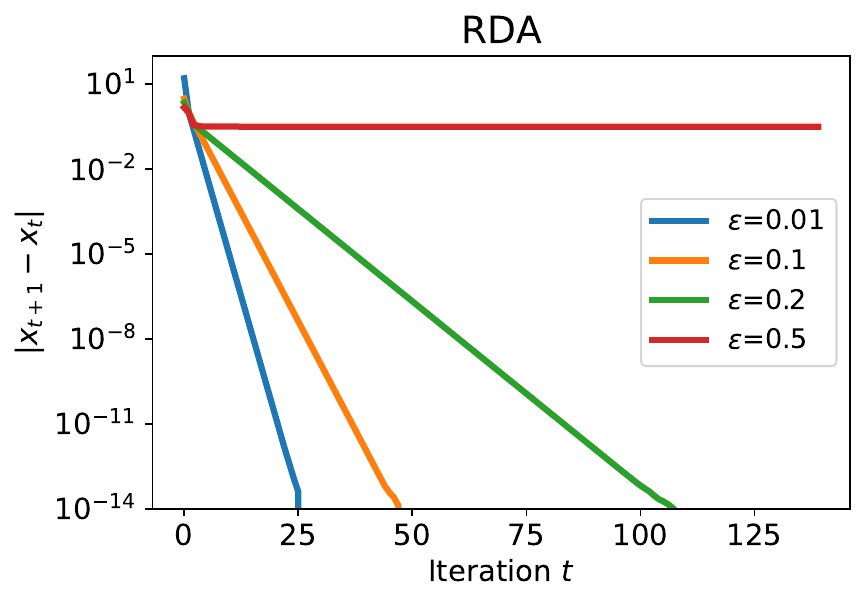}}
\caption{Convergence of RDA for different $\epsilon$-sensitivity parameters in the dynamic parking pricing problem.}
\label{RDA:pricing_simulation}
\end{center}
\vskip -0.2in
\end{figure}

We use the data from the street of Beach ST 600, containing a total of $n=10224$ samples.
Six features are considered, with the total occupied time treated as the strategic feature.
The occupancy is computed by dividing the total time by the total occupied time.
As in \cite{ray2022decision}, we approximate the distribution of occupancy as $z_1 = \zeta - A x$, where $A\approx 0.157$ is estimated from data,  $\zeta\sim p_0$ is a fixed distribution and $p_0$ is sampled from data of which the price is at the nominal price.

The constraint is selected as $\mathop{\mathbb{E}}\limits_{z_2 \sim \mathcal{D}(x)}[z_2] \leq c_1 x +c_2$, meaning that the total occupied time should be linearly bounded. We note that also the distribution in the constraint is $\epsilon$-sensitive.
%
At iteration $t$, the decision-dependent problem can be written as 
\begin{align*}
    \min_x & \;  (\frac{v}{2} +A^2) x^2  -2A \mathop{\mathbb{E}}_{\zeta\sim p_0} [\zeta]  x +1.4A x -t \mathop{\mathbb{E}[ z_2]}_{z_2\sim \mathcal{D}(x_t)} \\
     {\rm{s.t.}}& \; -c_1 x\leq  c_2- \mathop{\mathbb{E}}\limits_{z_2 \sim \mathcal{D}(x_{t})}[z_2].
\end{align*}

We select $v=0.03$, $t=0.005$, $c_1 =0.5$, $c_2 = 5$, $\eta = 0.18$ and $x_0 =-1$. All the minimization problems are solved using the CVXOPT toolbox. The convergence results for RCM and RDA are presented in Figs.~\ref{RCM:pricing_simulation} and \ref{RDA:pricing_simulation}, respectively.
In Fig.~\ref{RCM:pricing_simulation}, we observe that RCM does indeed converge in only a few iterations for
small values of $\epsilon$ while it divergences if $\epsilon$ is too large. Fig.~\ref{RDA:pricing_simulation} shows that
RDA also converges linearly for small values of $\epsilon$, however, with a slower rate compared to RCM.

\section{Conclusion}\label{sec:conclusion}
In this work, we developed a framework for solving constrained optimization with decision-dependent distributions in both the objective function and linear constraints.
Firstly, we established a sufficient condition that guarantees the uniqueness of the constrained equilibrium point and at the same time ensures that the RCM algorithm converges to this point.
Furthermore, we presented sufficient conditions for the RDA algorithm to converge to the constrained equilibrium point.
Additionally, we derived a bound for the distance between the constrained equilibrium and optimal points.
Our results showed that the decision-dependent optimization in \cite{perdomo2020performative} can be viewed as a special case of our framework.
Finally, we demonstrated the effectiveness of our algorithms using illustrative experiments on both a synthetic market problem and a dynamic pricing problem based on an open-source dataset.

\appendix\label{sec:appendix}
\subsection{Proof of Theorem~\ref{thm:RCM}}\label{app:RCM}
For ease of notation, we define
\begin{equation*}
    M(a,b) = \mathop{\rm{arg \; min}}_{x}  \mathop{\mathbb{E}}_{z\sim \mathcal{D}(a)}[l(x,z)] \quad {\rm{s.t.}} \quad  Gx \leq \xi(b).
\end{equation*}
We note that $M(x',x')$ and $M(x'',x')$ are minima of the objective  functions $\mathop{\mathbb{E}}_{z\sim \mathcal{D}(x')}[l(x,z)]$ and $\mathop{\mathbb{E}}_{z\sim \mathcal{D}(x'')}[l(x,z)]$, respectively, under the same constrained set. Hence, we can use Lemma \ref{lemma:PP} to bound the distance of these two minima, which gives
\begin{align}\label{eq:th1:p1}
    \left\| M(x',x')-M(x'',x')\right\| \leq \frac{\epsilon \beta_z}{\gamma} \left\| x' - x'' \right\|.
\end{align}

The two points $M(x'',x'')$ and $M(x'',x')$ are minima of the same objective function in two different constraint sets.
By virtue of Lemma~\ref{lemma:sensitivity}, the distance between $M(x'',x'')$ and $M(x'',x')$ can be bounded by 
\begin{align}\label{eq:th1:p2}
    \left\| M(x'',x'')-M(x'',x')\right\| &\leq L_{x^{*}} \left\| \xi(x') - \xi(x'') \right\| \nonumber \\
    &\leq L_{x^{*}} \epsilon_g \left\| x' - x'' \right\|.
\end{align}

Combining \eqref{eq:th1:p1} and \eqref{eq:th1:p2}, we have
\begin{align}\label{eq:RCM:contraction}
    &\left\| M(x',x')-M(x'',x'')\right\| \nonumber \\
    &\leq \left\| M(x',x')-M(x'',x')\right\| \nonumber \\
    &\quad + \left\| M(x'',x'')-M(x'',x')\right\| \nonumber \\
    &\leq \Big(  \frac{\epsilon \beta_z}{\gamma}+ L_{x^{*}} \epsilon_g\Big)\left\| x' - x'' \right\|.
\end{align}
From the update equation \eqref{eq:RCM}, we have $x_{t+1} = M(x_t,x_t)$ and $x_s = M(x_s,x_s)$. 
From \eqref{eq:RCM:contraction}, we have
\begin{align}
    \left\| x_t -x_s \right\| &= \left\| M(x_{t-1},x_{t-1}) -M(x_s,x_s) \right\| \nonumber \\
    &\leq \Big(  \frac{\epsilon \beta_z}{\gamma}+ L_{x^{*}} \epsilon_g\Big) \left\| x_{t-1} -x_s \right\| \nonumber \\
    &\leq \Big(  \frac{\epsilon \beta_z}{\gamma}+ L_{x^{*}} \epsilon_g\Big)^t \left\| x_0-x_s\right\| \nonumber ,
\end{align}
which shows that RCM converges at a linear rate if $\frac{\epsilon \beta_z}{\gamma}+ L_{x^{*}} \epsilon_g<1$.
Given that the mapping $x \rightarrow M(x,x)$ is a contraction, the uniqueness of $x_s$ directly follows from the Banach fixed-point theorem. The proof is complete.

\subsection{Proof of Theorem~\ref{thm:RDA}}\label{app:RDA}
We first note that all the items are positive since $\rho_2 (1+\rho_3) <2\gamma_d$ and $\rho_4 <1$.

By using the KKT conditions, it is easy to show that
\begin{align}
    x_s &= \mathop{\rm{arg \; min}}_{x\in \mathbb{R}^n} L_{x_s}(x,\lambda_s^{*})  \nonumber\\
    \lambda_s^{*}&= [\lambda_s^{*} + \eta \nabla d_{x_s}(\lambda_s^{*})]_{+},
\end{align}
where $\lambda_s^{*}$ is the optimal multiplier with respect to $d_{x_s}(\lambda)$. 
We first analyze the dynamics of the dual variable. From the update rule \eqref{eq:RDA:dual:update}, we have
\begin{align}\label{eq:th2:p1}
    &\left\|\lambda_{t+1} - \lambda_s^{*}\right\|^2 \nonumber \\
    &= \left\|[\lambda_t + \eta \nabla d_{x_t}(\lambda_t)]_{+} - [\lambda_s^{*} + \eta \nabla d_{x_s}(\lambda_s^{*})]_{+} \right\|^2 \nonumber \\
    &\leq  \left\|[\lambda_t + \eta \nabla d_{x_t}(\lambda_t)] - [\lambda_s^{*} + \eta \nabla d_{x_s}(\lambda_s^{*})] \right\|^2 \nonumber \\
    &=  \left\| \lambda_t - \lambda_s^{*} \right\|^2  + \eta^2 \left\|\nabla d_{x_t}(\lambda_t) - \nabla d_{x_s}(\lambda_s^{*}) \right\|^2 \nonumber \\
    &\quad + 2\eta \langle\lambda_t - \lambda_s^{*}, \nabla d_{x_t}(\lambda_t) - \nabla d_{x_s}(\lambda_s^{*}) \rangle.  
\end{align}
Since $d_{x'}(\lambda)$ is $\gamma_d$-strongly concave for every $x'$, we have 
\begin{align}\label{eq:th2:p2}
    & \langle\lambda_t - \lambda_s^{*}, \nabla d_{x_t}(\lambda_t) - \nabla d_{x_s}(\lambda_s^{*}) \rangle \nonumber \\
    &=  \langle\lambda_t - \lambda_s^{*}, \nabla d_{x_t}(\lambda_t) - \nabla d_{x_s}(\lambda_t) \rangle\nonumber \\
    &\quad + \langle\lambda_t - \lambda_s^{*},\nabla d_{x_s}(\lambda_t)- \nabla d_{x_s}(\lambda_s^{*}) \rangle \nonumber \\
    &\leq  \langle\lambda_t - \lambda_s^{*}, \nabla d_{x_t}(\lambda_t) - \nabla d_{x_s}(\lambda_t)\rangle - \gamma_d \left\| \lambda_{t} - \lambda_s^{*}\right\|^2 \nonumber \\
    &\leq  \left\| \lambda_{t} - \lambda_s^{*}\right\| \left\| \nabla d_{x_t}(\lambda_t) - \nabla d_{x_s}(\lambda_t)\right\| - \gamma_d \left\| \lambda_{t} - \lambda_s^{*}\right\|^2.
\end{align}
Moreover, we have
\begin{align}\label{eq:th2:p3}
    &\left\| \nabla d_{x_t}(\lambda_t) - \nabla d_{x_s}(\lambda_t)\right\| \nonumber \\
    & =  \Big\| G \nabla \tilde{f}_{x_t}(-G^{\rm{T}} \lambda_t) - \xi(x_t)  - G \nabla \tilde{f}_{x_s}(-G^{\rm{T}} \lambda_t) + \xi(x_s)\Big\| \nonumber \\
    & =  \left\| G \bar{y}_t - G \tilde{y}_t - \xi(x_t)   + \xi(x_s)\right\|,
\end{align}
where $ \bar{y}_t = \mathop{\rm{arg \; min}}\limits_{x \in \mathbb{R}^n } (f_{x_t}(x) + \lambda_t^{\rm{T}} Gx - \xi(x_t))$, $ \tilde{y}_t = \mathop{\rm{arg \; min}}\limits_{x \in \mathbb{R}^n } (f_{x_s}(x) + \lambda_t^{\rm{T}} Gx - \xi(x_s))$. By virtue of Lemma~\ref{lemma:PP}, we have 
\begin{align}\label{eq:th2:p4}
    \left\|\bar{y}_t -\tilde{y}_t\right\|\leq \frac{\epsilon \beta_z}{\gamma}\left\| x_t - x_s \right\| = \rho_1 \left\| x_t - x_s \right\|.
\end{align}
Substituting \eqref{eq:th2:p4} into \eqref{eq:th2:p3}, we obtain
\begin{align}\label{eq:th2:p5}
    &\left\| \nabla d_{x_t}(\lambda_t) - \nabla d_{x_s}(\lambda_t)\right\| \nonumber \\
    &\leq \left(\left\| G \right\|_2\frac{\epsilon \beta_z}{\gamma} + \epsilon_g \right)\left\| x_t - x_s\right\| = \rho_2 \left\| x_t - x_s\right\|.
\end{align}
Combining \eqref{eq:th2:p5} with \eqref{eq:th2:p2}, we have
\begin{align}\label{eq:th2:p6}
    &\langle\lambda_t - \lambda_s^{*}, \nabla d_{x_t}(\lambda_t) - \nabla d_{x_s}(\lambda_s^{*})\rangle \nonumber \\
    &\leq \rho_2 \left\| x_t - x_s\right\|  \left\| \lambda_{t} - \lambda_s^{*}\right\|- \gamma_d \left\| \lambda_{t} - \lambda_s^{*}\right\|^2.
\end{align}
Besides, we have
\begin{align}\label{eq:th2:p7}
    &\left\|\nabla d_{x_t}(\lambda_t) - \nabla d_{x_s}(\lambda_s^{*}) \right\|^2 \nonumber \\
    & =   \left\|\nabla d_{x_t}(\lambda_t) - \nabla d_{x_s}(\lambda_t) + \nabla d_{x_s}(\lambda_t) -\nabla d_{x_s}(\lambda_s^{*}) \right\|^2  \nonumber \\
    &= \left\|\nabla d_{x_t}(\lambda_t) - \nabla d_{x_s}(\lambda_t)\right\|^2 + \left\| \nabla d_{x_s}(\lambda_t) -\nabla d_{x_s}(\lambda_s^{*}) \right\|^2  \nonumber \\
    &\quad + 2 \langle\nabla d_{x_t}(\lambda_t) - \nabla d_{x_s}(\lambda_t),  \nabla d_{x_s}(\lambda_t) -\nabla d_{x_s}(\lambda_s^{*})\rangle  \nonumber \\
    &\leq   \rho_2^2 \left\| x_t - x_s\right\|^2 + \left\| \nabla d_{x_s}(\lambda_t) -\nabla d_{x_s}(\lambda_s^{*}) \right\|^2 \nonumber \\
    &\quad + 2  \left\| \nabla d_{x_t}(\lambda_t) - \nabla d_{x_s}(\lambda_t)\right\| \left\| \nabla d_{x_s}(\lambda_t) -\nabla d_{x_s}(\lambda_s^{*})\right\| \nonumber \\
    &\leq  \rho_2^2  \left\| x_t - x_s\right\|^2 + L_d^2 \left\| \lambda_{t} - \lambda_s^{*}\right\|^2 \nonumber \\
    &\quad + 2 L_d \rho_2 \left\| \lambda_{t} - \lambda_s^{*}\right\| \left\| x_t - x_s\right\|,
\end{align}
where the first inequality follows from the Cauchy–Schwarz inequality. The last inequality follows from the smoothness of $d_{x_s}(\lambda)$, which gives $\left\| \nabla d_{x_s}(\lambda_t) -\nabla d_{x_s}(\lambda_s^{*}) \right\|\leq L_d\left\| \lambda_{t} - \lambda_s^{*}\right\| $. From \eqref{eq:th2:p1}, \eqref{eq:th2:p6}, and \eqref{eq:th2:p7}, we have
\begin{align}\label{eq:th2:p8}
    &\left\|\lambda_{t+1} - \lambda_s^{*}\right\|^2  \nonumber \\
    &\leq   \left\| \lambda_t - \lambda_s^{*} \right\|^2  + \eta^2 \Big( \rho_2^2  \left\| x_t - x_s\right\|^2 + L_d^2 \left\| \lambda_{t} - \lambda_s^{*}\right\|^2 \nonumber \\
    &\quad \quad + 2 L_d \rho_2 \left\| \lambda_{t} - \lambda_s^{*}\right\| \left\| x_t - x_s\right\|\Big) \nonumber \\
    &\quad +2\eta \Big( \rho_2 \left\| x_t - x_s\right\|  \left\| \lambda_{t} - \lambda_s^{*}\right\|- \gamma_d \left\| \lambda_{t} - \lambda_s^{*}\right\|^2\Big) \nonumber \\
    &= \Big( 1 +\eta^2 L_d^2 - 2\eta \gamma_d \Big) \left\| \lambda_t - \lambda_s^{*} \right\|^2 +\eta^2 \rho_2^2  \left\| x_t - x_s\right\|^2 \nonumber \\
    &\quad +2\rho_2(L_d \eta^2 +\eta) \left\| x_t - x_s\right\|  \left\| \lambda_{t} - \lambda_s^{*}\right\| \nonumber \\
    &\leq  \Big( 1 +\eta^2 L_d^2 - 2\eta \gamma_d +\rho_2 \eta +\rho_2L_d \eta^2 \Big) \left\| \lambda_t - \lambda_s^{*} \right\|^2 \nonumber \\
    &\quad + \Big( \eta^2 \rho_2^2 +\rho_2 \eta +\rho_2 L_d \eta^2 \Big)\left\| x_t - x_s\right\|^2,
\end{align}
where the last inequality follows from the fact that $2ab \leq a^2 + b^2$.

Next we analyze the dynamics of the primal variable.
Recall that 
\begin{align*}
    x_t& =\mathop{\rm{arg \; min}}_{x} (f_{x_{t-1}}(x) + \lambda_t^{\rm{T}} (Gx - \xi(x_{t-1}))) \nonumber \\
    &=\mathop{\rm{arg \; min}}_{x} \mathop{\mathbb{E}}_{z\sim \mathcal{D}(x_{t-1})}[l(x,z) + \lambda_t^{\rm{T}}Gx ]\nonumber\\
    x_s &=\mathop{\rm{arg \; min}}_{x} (f_{x_{s}}(x) + \lambda_s^{*{\rm{T}}} (Gx - \xi(x_{s}))) \nonumber \\
    &=\mathop{\rm{arg \; min}}_{x} \mathop{\mathbb{E}}_{z\sim \mathcal{D}(x_{s})}[l(x,z) +\lambda_s^{*{\rm{T}}} Gx],
\end{align*}
and denote 
\begin{align*}
    y_t :&=\mathop{\rm{arg \; min}}_{x} (f_{x_{s}}(x) + \lambda_t^{\rm{T}} (Gx - \xi(x_{s})))\nonumber \\
    &=\mathop{\rm{arg \; min}}_{x} \mathop{\mathbb{E}}_{z\sim \mathcal{D}(x_{s})}[l(x,z) +\lambda_t^{\rm{T}} Gx]
\end{align*}
It can be verified that the function $l(x,z)+ \lambda_t^{\rm{T}} Gx$ is $\gamma$-strongly convex in $x$ and its gradient $\nabla_x l(x,z) + G^T \lambda_t$ is $\beta$-Lipschitz continuous in $z$.
By virtue of Lemma \ref{lemma:PP}, the distance between $x_t$ and $y_t$ can be bounded by
\begin{align}\label{eq:th2:p9}
    \left\| x_{t} -y_t\right\| \leq \frac{\epsilon \beta_z}{\gamma} \left\| x_{t-1}-x_s\right\| = \rho_1 \left\| x_{t-1}-x_s\right\|.
\end{align}
Now we bound $\left\| y_t -x_s\right\|$. 
Since $f_{x_s}(x)$ is strongly convex in $x$, we have
\begin{align}\label{eq:th2:p10}
    &\gamma \left\| x_s - y_t\right\|^2 \leq \langle x_s-y_t, \nabla f_{x_s}(x_s) - \nabla f_{x_s}(y_t)\rangle \nonumber \\
    &\leq \left\| x_s - y_t \right\|  \left\| \nabla f_{x_s}(x_s) - \nabla f_{x_s}(y_t) \right\|.
\end{align}
Since $f_{x_{s}}(x) + \lambda_t^{\rm{T}} Gx$ is also strongly convex in $x$, $y_t$ is unique and satisfies $\nabla f_{x_s}(y_t) + G^{\rm{T}} \lambda_t = 0$. Similarly, $f_{x_{s}}(x) + \lambda_s^{*{\rm{T}}} Gx$ is strongly convex in $x$ and thus $\nabla f_{x_{s}}(x_s) + G^{\rm{T}} \lambda_s^{*}=0$. 
Then, \eqref{eq:th2:p10} yields
\begin{align}\label{eq:th2:p11}
    \left\| x_s - y_t\right\| &\leq \frac{1}{\gamma} \left\| \nabla f_{x_s}(x_s) - \nabla f_{x_s}(y_t) \right\|  \nonumber \\
    &= \frac{1}{\gamma} \left\| -G^{\rm{T}}\lambda_s^{*}+ G^{\rm{T}} \lambda_t\right\| \nonumber \\
    &\leq \frac{\left\| G \right\|_2}{\gamma}  \left\| \lambda_t - \lambda_s^{*} \right\|.
\end{align}
From \eqref{eq:th2:p9} and \eqref{eq:th2:p11}, we have
\begin{align}\label{eq:th2:p12}
    & \left\| x_t - x_s \right\|^2 = \left\| x_t - y_t +y_t - x_s \right\|^2 \nonumber \\
    &=   \left\| x_t - y_t \right\|^2 + \left\| y_t - x_s \right\|^2 + 2\langle x_t-y_t,y_t-x_s\rangle \nonumber \\
    &\leq  \rho_1^2\left\| x_{t-1}-x_s\right\|^2  +\frac{\left\| G \right\|_2^2}{\gamma}  \left\| \lambda_t - \lambda_s^{*} \right\|^2 \nonumber \\
    &\quad + 2  \rho \frac{\left\| G \right\|_2}{\gamma}\left\| x_{t-1}-x_s\right\|\left\| \lambda_t - \lambda_s^{*} \right\| \nonumber \\
    &\leq  \Big(\rho_1^2 + \frac{\rho_1\left\| G \right\|_2 }{\gamma}\Big) \left\| x_{t-1}-x_s\right\|^2 \nonumber \\
    &\quad + \Big( \frac{\left\| G \right\|_2^2}{\gamma^2}  +\frac{\rho_1\left\| G \right\|_2 }{\gamma} \Big) \left\| \lambda_t - \lambda_s^{*} \right\|^2 \nonumber \\
    &=\rho_4 \left\| x_{t-1}-x_s\right\|^2 + \rho_3 \left\| \lambda_t - \lambda_s^{*} \right\|^2.
\end{align}
Combining \eqref{eq:th2:p8} and \eqref{eq:th2:p12}, we have 
\begin{align}\label{eq:th2:p13}
    &\left\|\lambda_{t+1} - \lambda_s^{*}\right\|^2 + \alpha  \left\| x_t - x_s \right\|^2 \nonumber \\
    &\leq  \Big( 1 +\eta^2 L_d^2 - 2\eta \gamma_d +\rho_2 \eta +\rho_2L_d \eta^2 \Big) \left\| \lambda_t - \lambda_s^{*} \right\|^2 \nonumber \\
    &\quad + \Big( \eta^2 \rho_2^2 +\rho_2 \eta +\rho_2 L_d \eta^2 
    + \alpha \Big)\left\| x_t - x_s\right\|^2 \nonumber \\
    &\leq \Big( 1 +\eta^2 L_d^2 - 2\eta \gamma_d +\rho_2 \eta +\rho_2L_d \eta^2 \Big) \left\| \lambda_t - \lambda_s^{*} \right\|^2 \nonumber \\
    &\quad + \Big( \eta^2 \rho_2^2 +\rho_2 \eta +\rho_2 L_d \eta^2 
    + \alpha \Big) \\
    & \quad \quad \times \Big( \rho_4 \left\| x_{t-1}-x_s\right\|^2 + \rho_3 \left\| \lambda_t - \lambda_s^{*} \right\|^2\Big) \nonumber \\
   & =  \Big( 1 +\eta^2 L_d^2 - 2\eta \gamma_d +\rho_2 \eta +\rho_2L_d \eta^2 \nonumber \\
    &\quad + \rho_3 \Big( \eta^2 \rho_2^2 +\rho_2 \eta +\rho_2 L_d \eta^2 
    + \alpha \Big)\Big) \left\| \lambda_t - \lambda_s^{*} \right\|^2 \nonumber \\
    &\quad + \frac{\rho_4}{\alpha} \Big( \eta^2 \rho_2^2 +\rho_2 \eta +\rho_2 L_d \eta^2 
    + \alpha \Big) \alpha\left\| x_{t-1}-x_s\right\|^2 \nonumber \\
    &= (1+H_{\lambda}(\eta) )\left\| \lambda_t - \lambda_s^{*} \right\|^2 \nonumber \\
    &\quad + (1+H'_{x}(\eta)) \alpha\left\| x_{t-1}-x_s\right\|^2,
\end{align}
where we define $H_{\lambda}(\eta) := \eta^2 L_d^2 - 2\eta \gamma_d +\rho_2 \eta +\rho_2L_d \eta^2 + \rho_3 \Big( \eta^2 \rho_2^2 +\rho_2 \eta +\rho_2 L_d \eta^2  + \alpha \Big)  = a_1 \eta^2 - a_2 \eta + a_3$ and $H'_{x}(\eta)  =  \frac{\rho_4}{\alpha} \Big( \eta^2 \rho_2^2 +\rho_2 \eta +\rho_2 L_d \eta^2 + \alpha \Big) -1 = \frac{1}{\alpha} \Big( b_1 \eta^2 +b_2 \eta - b_3 \Big) := \frac{1}{\alpha} H_{x}(\eta) $. It suffices to show that the choice of $\eta\in(s_1,s_2)$ assures $H_{\lambda}(\eta) <0$ and  $H_{x}(\eta)<0$.

Since $a_2^2 > 4a_1a_2$, there exist two distinct solutions to $H_{\lambda}(\eta) =0$, which we denote by $s_1$, $s'_1$ with $s_1 <s'_1$. When $\eta \in (s_1,s'_1)$, we have $H_{\lambda}(\eta) <0$. Similarly, we denote by $s_2$, $s'_2$ the solutions to $H_{x}(\eta)<0$ with $s_2>s'_2$. When $\eta \in (s'_2,s_2)$, we have $H_{\lambda}(x) <0$. When $\eta$ takes values on the intersection of these two regions, we can make sure that $H_{\lambda}(\eta) <0$ and $H_{x}(\eta)<0$ simultaneously. Obviously, this happens if $s_1<s_2$. In what follows, we show that the choice of $\alpha$ ensures $s_1<s_2$.
Specifically, the condition $s_1<s_2$ is equivalent to 
\begin{align}\label{eq:temp_21}
    &\frac{a_2 - \sqrt{a_2^2 - 4a_1 a_3}}{2a_1} < \frac{-b_2 + \sqrt{b_2^2 +4 b_1 b_3}}{2b_1} \nonumber \\
    \Leftrightarrow \qquad & a_2 b_1 + a_1 b_2 < b_1\sqrt{a_2^2 - 4a_1 a_3} + a_1\sqrt{b_2^2 +4 b_1 b_3} \nonumber \\
    \Leftrightarrow \qquad & \rho_2 \rho_4( L_d^2 + L_d \rho_2 + \rho_2^2 \rho_3 + L_d \rho_2 \rho_3) \nonumber \\
    &\quad + \rho_4(\rho_2^2 + L_d \rho_2) ( 2\gamma_d -\rho_2(1+\rho_3)) \nonumber \\
    &< b_1\sqrt{a_2^2 - 4a_1 a_3} + a_1\sqrt{b_2^2 +4 b_1 b_3} \nonumber \\
    \Leftrightarrow \qquad & \rho_2 \rho_4 (2\gamma_d \rho_2 + 2\gamma_d L_d + L_d^2 - \rho_2^2) \nonumber \\
    & <b_1\sqrt{a_2^2 - 4a_1 a_3} + a_1\sqrt{b_2^2 +4 b_1 b_3}.
\end{align}
Since $\sqrt{b_2^2 +4 b_1 b_3}>b_2$, a sufficient condition for \eqref{eq:temp_21} is 
\begin{align}\label{eq:temp_22}
    &\rho_2 \rho_4 (2\gamma_d \rho_2 + 2\gamma_d L_d + L_d^2 - \rho_2^2) \nonumber \\
    &<b_1\sqrt{a_2^2 - 4a_1 a_3} + a_1 b_2 \nonumber \\
    \Leftrightarrow \qquad &\rho_2 \rho_4 (2\gamma_d \rho_2 + 2\gamma_d L_d + L_d^2 - \rho_2^2) \nonumber \\
    &<b_1\sqrt{a_2^2 - 4a_1 a_3} + \rho_2 \rho_4(L_d^2 + L_d \rho_2 + \rho_2^2 \rho_3 \nonumber \\
    &\quad + L_d \rho_2 \rho_3) \nonumber \\
    \Leftrightarrow \qquad & 2\gamma_d\rho_2 + 2 \gamma_d L_d < (\rho_2 +L_d) \sqrt{a_2^2 - 4a_1 a_3} \nonumber \\
    &+ L_d^2 +L_d \rho_2 +\rho_2^3 \rho_3 + L_d\rho_2^2 \rho_3 + \rho_2^2 .
\end{align}
Since $\alpha = (1-\bar{\alpha})\frac{a_2^2}{4\rho_3 a_1}$, the condition \eqref{eq:temp_22} is equivalent to 
\begin{align}\label{eq:temp_23}
    & 2\gamma_d\rho_2 + 2 \gamma_d L_d < \sqrt{\bar{\alpha}}(\rho_2 +L_d) (2\gamma_d - \rho_2 (1+\rho_3)) \nonumber \\
    &+ L_d^2 +L_d \rho_2 +\rho_2^3 \rho_3 + L_d\rho_2^2 \rho_3 + \rho_2^2 \nonumber \\
    \Leftrightarrow \qquad &  2\gamma_d\rho_2 + 2 \gamma_d L_d + \sqrt{\bar{\alpha}} \rho_2 (1+\rho_3)(\rho_2 +L_d) \nonumber \\
    &<  2\sqrt{\bar{\alpha}}\gamma_d(\rho_2 +L_d) + L_d^2 +L_d \rho_2 +\rho_2^3 \rho_3 \nonumber \\
    &\quad + L_d\rho_2^2 \rho_3 + \rho_2^2 .
\end{align}
Since $\sqrt{\bar{\alpha}}\in [0,1)$ and $\gamma_d < L_d$, we further simplify the sufficient condition \eqref{eq:temp_23} to
\begin{align}
    &2\gamma_d \rho_2 + 2 \gamma_d L_d< \sqrt{\bar{\alpha}}(2\gamma_d \rho_2 + 2 \gamma_d L_d) + L_d^2\nonumber \\
    \Leftrightarrow \qquad & (1- \sqrt{\bar{\alpha}})(2\gamma_d \rho_2 + 2 \gamma_d L_d)<L_d^2 \nonumber \\
    \Leftrightarrow \qquad &  \sqrt{\bar{\alpha}} >1-\frac{L_d^2}{2\gamma_d(\rho_2 +L_d)},
\end{align}
which holds since $\sqrt{\bar{\alpha}} = \max\left\{0,1-\frac{L_d^2}{4\gamma_d(\rho_2 +L_d)}\right\}$.
Therefore, when $\eta\in(s_1,s_2)$, we have $H_x(\lambda):= \kappa_1<1$ and $H_x(\eta):= \kappa_2<1$. For $\kappa = \max\{\kappa_1,\kappa_2  \}$, from \eqref{eq:th2:p13}, we have $\left\| \lambda_{t+1} - \lambda_s^{*} \right\|^2 + \alpha \left\| x_t - x_s \right\|^2 \leq \kappa (\left\| \lambda_{t} - \lambda_s^{*} \right\|^2 + \alpha \left\| x_{t-1} - x_s \right\|^2) $. Iteratively using this inequality completes the proof.

\subsection{Proof of Corollary~\ref{corollary:RDA}}\label{app:corollary2}
Following the steps as in the proof of Theorem~\ref{thm:RDA}, we have that
\begin{align}\label{eq:corollary2:t1}
    &\left\|\lambda_{t+1} - \lambda_s^{*}\right\|^2  \nonumber \\
    &\leq  \Big( 1 +\eta^2 L_d^2 - 2\eta \gamma_d +\rho_2 \eta +\rho_2L_d \eta^2 \Big) \left\| \lambda_t - \lambda_s^{*} \right\|^2 \nonumber \\
    &\quad + \Big( \eta^2 \rho_2^2 +\rho_2 \eta +\rho_2 L_d \eta^2 \Big)\left\| x_t - x_s\right\|^2,
\end{align}
and 
\begin{align}\label{eq:corollary2:t2}
    \left\| x_{t+1} - x_s \right\|^2 \leq \frac{\left\| G\right\|_2^2}{\gamma^2} \left\| \lambda_t - \lambda_s^{*}\right\|^2.
\end{align}
Substituting \eqref{eq:corollary2:t2} into \eqref{eq:corollary2:t1}, we have
\begin{align}\label{eq:corollary2:t3}
    &\left\| \lambda_{t+1}- \lambda_s^{*} \right\|^2 \nonumber \\
    &\leq \Big( 1 - 2\eta \gamma_d +\eta^2 L_d^2 + \epsilon_g \eta + \epsilon_g L_d \eta^2  \nonumber \\
    &\quad \quad  +\frac{\left\| G\right\|_2^2 } {\gamma^2} (\epsilon_g \eta^2 + \epsilon_g \eta + \epsilon_g L_d \eta^2) \Big) \left\| \lambda_t - \lambda_s^{*}\right\|^2 \nonumber \\
    & = \Big( 1 - \eta(2\gamma_d - \epsilon_g - \epsilon_g \frac{\left\| G\right\|_2^2}{\gamma^2} )  +\eta^2 (L_d^2 + \epsilon_g L_d \nonumber \\
    &\quad \quad  + \epsilon_g  \frac{\left\| G\right\|_2^2}{\gamma^2} + \epsilon_g L_d \frac{\left\| G\right\|_2^2}{\gamma^2}) \Big) \left\| \lambda_t - \lambda_s^{*}\right\|^2 \nonumber \\
    &= \kappa_3 \left\| \lambda_t - \lambda_s^{*}\right\|^2.
\end{align}
Since $\eta <\frac{2\gamma_d - \epsilon_g (1 + \frac{\left\| G \right\|_2^2}{\gamma^2} )}{L_d^2 + \epsilon_g L_d + \epsilon_g^2 \frac{\left\| G \right\|_2^2}{\gamma^2} + \epsilon_g^2 L_d \frac{\left\| G \right\|_2^2}{\gamma^2} }$, we have $\kappa_3 < 1$. Iteratively using \eqref{eq:corollary2:t3}, we can obtain \eqref{eq:corollary:lambda_converge}. 
Substituting \eqref{eq:corollary2:t2} into \eqref{eq:corollary:lambda_converge}, we obtain \eqref{eq:corollary:x_converge}. The proof is complete.

\subsection{Proof of Theorem~\ref{thm:optimality_stability}}\label{app:optimality_stability}
Recall the definition $f_{x'}(x) = \mathop{\mathbb{E}}\limits_{z\sim \mathcal{D}(x')}[l(x,z)]$.
By the definition of $x_o$ and $x_s$, we have $G x_o \leq \xi(x_o)$ and $G x_s \leq \xi(x_s)$. Since both $x_o$ and $x_s$ are feasible points of the problem \eqref{eq:problem}, we have $f_{x_o}(x_o)\leq f_{x_s}(x_s)$.
We note that the inequality $ f_{x_s}(x_s)\leq f_{x_s}(x_o)$ does not necessarily hold since
the constraints $G x_o \leq \xi(x_s)$ do not necessarily hold. 
In fact, we can only guarantee that 
\begin{align}\label{eq:thm3:pf:1}
    G_i x_o &\leq \xi_i(x_o) -\xi_i(x_s) +\xi_i(x_s) \nonumber \\
    &\leq \xi_i(x_s) + \epsilon_g \left\| x_o-x_s\right\|,
\end{align}
for all $i\in \{ 1,\ldots, d_w\}$, where $\xi_i(x_s)$ denotes the $i$-th element of the vector $\xi(x_s)$.

To bound $\left\| x_o - x_s \right\| $, we define the point 
\begin{align*}
    \hat{x}_s =  &\mathop{\rm{arg \; min}}_{x} \Big\{ \mathop{\mathbb{E}}_{z\sim \mathcal{D}(x_s)}[l(x,z)] \Big| \nonumber \\
    &G_i x\leq \xi_i(x_s) +\epsilon_g \left\| x_o-x_s\right\|, i\in \{1,\ldots,d_w\} \Big\}.
\end{align*}
Since both $\hat{x}_s$ and $x_o$ satisfy the constraint \eqref{eq:thm3:pf:1}, for all $i\in \{ 1,\ldots, d_w\}$, the strongly convex property of $ f_{x'}(x)$ yields
\begin{align}\label{eq:th3:p1}
    f_{x_s}(x_o) - f_{x_s}(\hat{x}_s) \geq \frac{\gamma}{2} \left\| x_o - \hat{x}_s \right\|^2.
\end{align}
Now we bound $|f_{x_s}(x_s) - f_{x_s}(\hat{x}_s)|$. Denote by $\lambda^{*}_s$ the optimal multiplier of the problem $\max\limits_{\lambda \geq 0} g_{x_s}(\lambda)$, where 
$$g_{x_s}(\lambda) = \min\limits_x \left\{\mathop{\mathbb{E}}\limits_{z\sim \mathcal{D}(x_s)}[l(x,z)] + \lambda^T(Gx - \xi(x_s))\right\}.$$ 
Consider the convex minimization problem:
\begin{align*}
    \min_{x} &\mathop{\mathbb{E}}\limits_{z\sim \mathcal{D}(x_s)}[l(x,z)] \nonumber \\
    \text{s.t.} &\quad  G_ix\leq \xi_i(x_s)+u, \; i=1,\ldots, d_w,
\end{align*}
where $u$ is a scalar.
Define $p^{*}(u) = \inf\limits_{x} \Big\{ \mathop{\mathbb{E}}\limits_{z\sim \mathcal{D}(x_s)}[l(x,z)] \Big|   G_i x\leq \xi_i(x_s)+u , i\in \{ 1,\ldots, d_w\}\Big\}$.
Set $u_0 = \epsilon_g \left\| x_o-x_s\right\|$. For any feasible point $x$ that satisfies the constraint $G_i x\leq \xi(x_s)+u_0$, $i=1,\ldots,d_w$, we have 
\begin{align}
    &f_{x_s}(x_s) = p^{*}(0) = g_{x_s}(\lambda_s^{*}) \nonumber \\
    &\leq \mathop{\mathbb{E}}\limits_{z\sim \mathcal{D}(x_s)}[l(x,z)] + \lambda_s^{*{\rm{T}}}(Gx - \xi(x_s))  \nonumber \\
    & \leq \mathop{\mathbb{E}}\limits_{z\sim \mathcal{D}(x_s)}[l(x,z)] + \lambda_s^{*{\rm{T}}} \textbf{1} u_0, \nonumber 
\end{align}
where the equality holds due to the strong duality, the first inequality follows from the definition of $g_{x_s}$, and the second inequality follows from the fact that $\lambda_s^{*}\geq 0$. Here, the notation $\textbf{1}$ denotes the vector with suitable dimension of which all the elements are 1.
Since $\hat{x}_s$ is a feasible point with respect to the constraint $G_i x\leq \xi_i(x_s)+u_0$, $i=1,\ldots,d_w$, we further obtain
\begin{align}\label{eq:th3:p2}
    f_{x_s}(x_s) &\leq \mathop{\mathbb{E}}\limits_{z\sim \mathcal{D}(x_s)}[l(\hat{x}_s,z)] + \lambda_s^{*{\rm{T}}} \textbf{1} u_0 \nonumber \\
    &\leq  f_{x_s}(\hat{x}_s) + \lambda_s^{*{\rm{T}}} \textbf{1} u_0 \nonumber \\
    &\leq f_{x_s}(\hat{x}_s)+ \epsilon_g \sqrt{d_w} \left\|\lambda_s^{*}\right\| \left\| x_o-x_s\right\|.
\end{align}
From \eqref{eq:th3:p1} and \eqref{eq:th3:p2}, we have
\begin{align}\label{eq:th3:p3}
    f_{x_o}(x_o) &\leq f_{x_s}(x_s) \nonumber \\
    &\leq f_{x_s}(\hat{x}_s) + \epsilon_g \sqrt{d_w}\left\| \lambda_s^{*} \right\|  \left\| x_o-x_s\right\| \nonumber \\
    &\leq  f_{x_s}(x_o) - \frac{\gamma}{2} \left\| x_o - \hat{x}_s \right\|^2 \nonumber \\
    &\quad + \epsilon_g \sqrt{d_w}\left\| \lambda_s^{*} \right\|  \left\| x_o-x_s\right\|.
\end{align}
Since $l(x,z)$ is $L_z$-Lipschitz continuous in $z$, we have
\begin{align}\label{eq:th3:p4}
    f_{x_s}(x_o)\leq f_{x_o}(x_o) +  L_z \epsilon \left\| x_o-x_s\right\|.
\end{align}
From \eqref{eq:th3:p3} and \eqref{eq:th3:p4}, we have
\begin{align}\label{eq:th3:p5}
    \frac{\gamma}{2} \left\| x_o - \hat{x}_s \right\|^2 
    \leq L_z \epsilon \left\| x_o-x_s\right\| \nonumber \\
    + \epsilon_g \sqrt{d_w}\left\| \lambda_s^{*} \right\|  \left\| x_o-x_s\right\|.
\end{align}
By virtue of Lemma \ref{lemma:sensitivity}, we have 
\begin{align}\label{eq:th3:p6}
    \left\| x_s - \hat{x}_s\right\| \leq L_{x^{*}} \epsilon_g \sqrt{d_w} \left\| x_o-x_s\right\|.
\end{align}
Combining \eqref{eq:th3:p5} with \eqref{eq:th3:p6}, we have
\begin{align}\label{eq:th3:p7}
    &\gamma \left\| x_o - x_s \right\|^2  = \gamma \left\|x_o - \hat{x}_s + \hat{x}_s- x_s  \right\|^2 \nonumber \\
    &\leq  \gamma (1+\frac{1}{c})\left\| x_o - \hat{x}_s\right\|^2 + \gamma (1+c) \left\| \hat{x}_s- x_s  \right\|^2 \nonumber \\
    & \leq   2(1+\frac{1}{c}) (L_z \epsilon + \epsilon_g \sqrt{d_w}\left\| \lambda_s^{*} \right\|)\left\| x_o-x_s\right\| \nonumber \\
    &\quad + \gamma (1+c) L_{x^{*}}^2 \epsilon_g^2 d_w \left\| x_o-x_s\right\|^2,
\end{align}
where the first inequality holds, since $2ab\leq c a^2 + \frac{1}{c} b^2$ for any $c>0$.
Upon performing a transformation, \eqref{eq:th3:p7} can be expressed as:
\begin{align*}
    \left\| x_o - x_s \right\| \leq  \frac{  2(1+\frac{1}{c}) (L_z \epsilon + \epsilon_g \sqrt{d_w} \left\| \lambda_s^{*} \right\|   }{ \gamma(1-(1+c) L_{x^{*}}^2 \epsilon_g^2 d_w)} .
\end{align*}
By setting $c=-1+\frac{1}{L_{x^{*}} \epsilon_g \sqrt{d_w}}$, we get 
\begin{align*}
     \left\| x_o - x_s \right\| \leq \frac{2(L_z \epsilon + \epsilon_g\sqrt{d_w} \left\| \lambda_s^{*} \right\|)}{\gamma (1-L_{x^{*}} \epsilon_g \sqrt{d_w})^2},
\end{align*}
which completes the proof.


\ifCLASSOPTIONcaptionsoff
  \newpage
\fi



%



\bibliography{00_citation}
\bibliographystyle{unsrt}








\end{document}